\documentclass[12pt,a4paper]{amsart}
\usepackage[T1]{fontenc}
\usepackage[utf8]{inputenc}
\usepackage{amsmath}
\usepackage{amsfonts}
\usepackage{amssymb}
\usepackage{setspace}
\usepackage{amsthm}
\usepackage{rotating}
\usepackage[a4paper, left=2.5cm, right=2.5cm, top=2.8cm, bottom=2.8cm]{geometry}
\usepackage{float}
\usepackage{amsmath}
\usepackage{arydshln}
\usepackage{graphicx}
\usepackage[usenames, dvipsnames]{color}
 \usepackage{tabu}
\usepackage{subfigure,bm}
\usepackage{natbib}
\usepackage{commath}
\usepackage{amssymb}
\usepackage{ mathrsfs }
\usepackage{booktabs} 
\usepackage{url} 
\usepackage{hyperref}
\usepackage{amsaddr} 

\newcommand\Ex{\text{\sf E}}
\newcommand\Prob{\text{\sf P}}

\newcommand\Cov{\text{\sf Cov}}

\newtheorem{theorem}{Theorem}
\newtheorem{proposition}{Proposition}
\newtheorem{lemma}{Lemma}
\newtheorem{corollary}{Corollary}

\theoremstyle{remark}

\newtheorem{remark}{Remark}

\begin{document}
\title[IIA through Slepian model]{Slepian model based independent interval approximation for the level excursion distributions}

\author[H.Bengtsson and K. Podg\'orski]{Henrik Bengtsson and Krzysztof Podg\'orski} \address{Department of Statistics,
Lund University, \\}
\email{Henrik.Bengtsson@stat.lu.se, \href{https://orcid.org/0000-0002-9280-4430}{ \includegraphics[height=2.2mm]{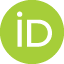} 0000-0002-9280-4430}}  
\email{Krzysztof.Podgorski@stat.lu.se, \href{https://orcid.org/0000-0003-0043-1532}{ \includegraphics[height=2.2mm]{ORCIDiD64.png} 0000-0003-0043-1532}}  
\date{January 27, 2025}

\begin{abstract}
The independent interval approximation of the excursion time distributions for Gaussian processes has been used in physics and engineering. A new but related approach matches the expected value of the clipped Slepian to the expected value of a non-stationary binary stochastic process. This approach is extended to non-zero crossings and provides a probabilistic foundation for the validity of the approximations for a large class of processes. Both the above and below distributions are approximated. While the Slepian-based method was shown to be equivalent to the classical IIA for the zero-level, this is not the case for non-zero excursions.
\end{abstract}

\maketitle

\section{Introduction}
\noindent
Excursion distribution has received considerable attention from the mid-1930s until today. This is not surprising, considering the importance of characterizing the distribution of excursions in many applications. For example, the length of a drought or a heatwave are two problems where characterizing the exceedance distributions is important. Here, classical methods in extreme value theory, such as peaks over threshold or block maxima, fail since it is the length of the excursion and not the maxima that is of interest.

The formulation of the excursion distribution is both simple and intuitive, but this hides the true complexity of the problem. The problem of fully characterizing the excursion distribution from characteristics of the underlying process is still unsolved despite considerable efforts. The study of excursion distributions originates in the study of the number of real roots of a random polynomial. It was perhaps \cite{LittlewoodOfford} who first formulated the question mathematically, and early work was done by \cite{kac1943}, {\cite{littlewood2}}, and \cite{Paul}. From this early work \cite{Rice} studied the expected number of roots of a stochastic process, which led to the famous Rice formula on the expected number of crossing in a unit of time for an ergodic stationary process. The study of the number of zeros would evolve into investigating the length between zeros, and hence, the first investigation of excursion times began. A discussion of the development is incomplete without mentioning the results of \cite{McFadden1956}.

\cite{McFadden1956} investigated the length of axis crossings and first proposed the notion of clipping a stochastic process. If $X(t)$ is the process of interest, the clipped process $Z(t)$ can be constructed by computing $Z(t)=sgn(X(t)-u)$. While $Z(t)$ will be less rich in terms of information, all the information regarding the excursion behavior of $X(t)$ will be contained in $Z(t)$.

Due to the lack of results connecting the characteristics of the process of interest to the excursion distribution, approximation methods have been developed. A method commonly used in statistical physics is the so-called independent interval approximation (IIA). The basis for this is to approximate the clipped process with a stationary binary process taking values one and minus one, where the time spent in the one and minus one states are independent. The extensive survey paper by \cite{BrayMS} is a good overview of the use of the IIA in statistical physics. However, \cite{LonguetHiggins1962} showed that this approximation method is never valid for Gaussian processes. Nevertheless, the IIA has been further developed, and in particular \cite{Sire2007} extended the IIA framework to non-zero levels. There are, however, questions about the mathematical validity of the IIA framework due to its lack of criterion on the characteristics used for the approximation.

The IIA has also been used to derive the so-called persistency coefficient. This coefficient describes the tail behavior of the excursion distribution and is used to characterize various systems in statistical physics. However, deriving the persistency coefficient is a difficult problem, and therefore, they are only known for a handful of processes. For example, \cite{PoplavskyiSchehr2018} derived the persistency coefficient for zero-level crossings of the two-dimensional diffusion process.

Recently, another version of the IIA framework was proposed by \cite{bengtsson2024slepian}. This version is based on clipping the Slepian process and matching it to a non-stationary switch process. However, this was only done for zero-level crossings, and this approach was shown to be equivalent to the ordinary IIA. In this paper, the Slepian-based IIA framework is extended to non-zero-level excursions. Conditions are also provided for the mathematical validity of this approach. These conditions are based on recent results on the properties of the switch process by \cite{BengtssonP}.

The paper is organized as follows. Section~\ref{sec:Slepian} and \ref{sec:switch} introduce the Slepian model and the Switch processes and thus serve as preliminaries for subsequent sections. In Section~\ref{sec:clipp}, properties of the clipped processes are derived, which are then used when matching characteristics. The latter is done in Section~\ref{sec:IIA}, and in this section, we demonstrate the difficulty of using the Clipped process and stationary switch process for non-zero level IIA. In Section~\ref{sec:diff}, the Slepian-based IIA is applied to approximate the excursion distribution of Gaussian diffusion processes in two dimensions. Persistency coefficients are approximated for this process, and the results are compared to results based on simulating trajectories of the process.

\section{Slepian process}
\label{sec:Slepian}
\noindent
Recall that the Slepian model at level $u$ corresponding to a smooth stationary Gaussian process $X(t)$ with the twice continuously differentiable normalized covariance function $r(t)$ is given through 
\begin{align}
\label{eq:Slepian}
X_{u}(t)=u \cdot r(t)-R \cdot {r'(t)}/{\sqrt{-r''(0)}}+\Delta(t),
\end{align}
where $r(t)$ is the covariance of $X$, $R$ is a standard Rayleigh variable, i.e. with density $f_R(s)=se^{-s^2/2}$ independent of a (non-stationary) Gaussian process $\Delta$ with covariance 
\begin{align*}
r_\Delta(t,s) = r(t-s)-r(t)r(s)+r'(t)r'(s)/r''(0).
\end{align*}
This process describes the statistical behavior of the original process $X$ at the instants of an up-crossing of the level $u$, i.e., the process behaves how a trajectory of a stationary Gaussian process is seen when observed at the $u$-level up-crossing instant.

\begin{figure}[t]
\includegraphics[width=0.43\textwidth]{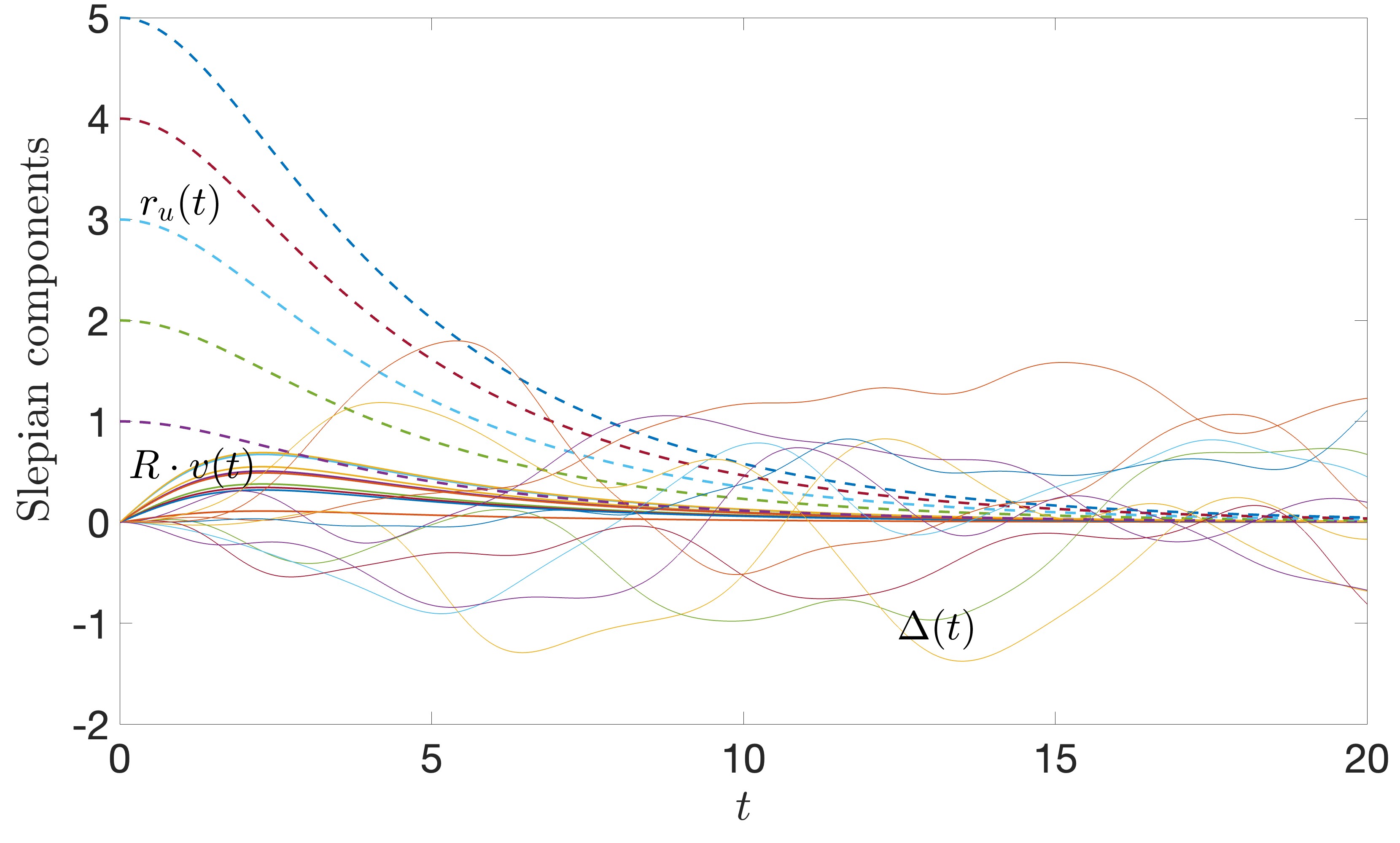}
\includegraphics[width=0.5\textwidth]{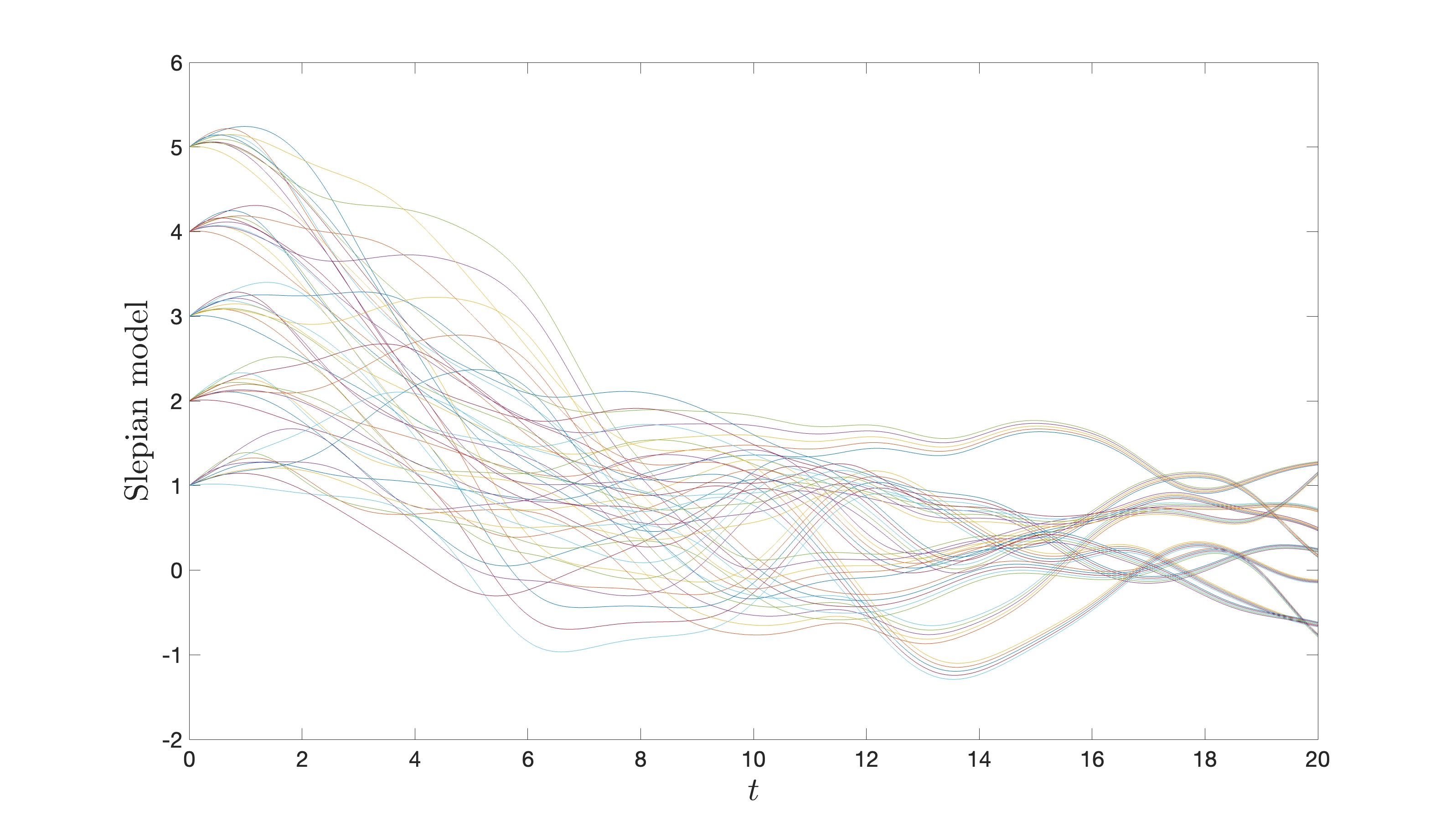}\\
\includegraphics[width=0.45\textwidth]{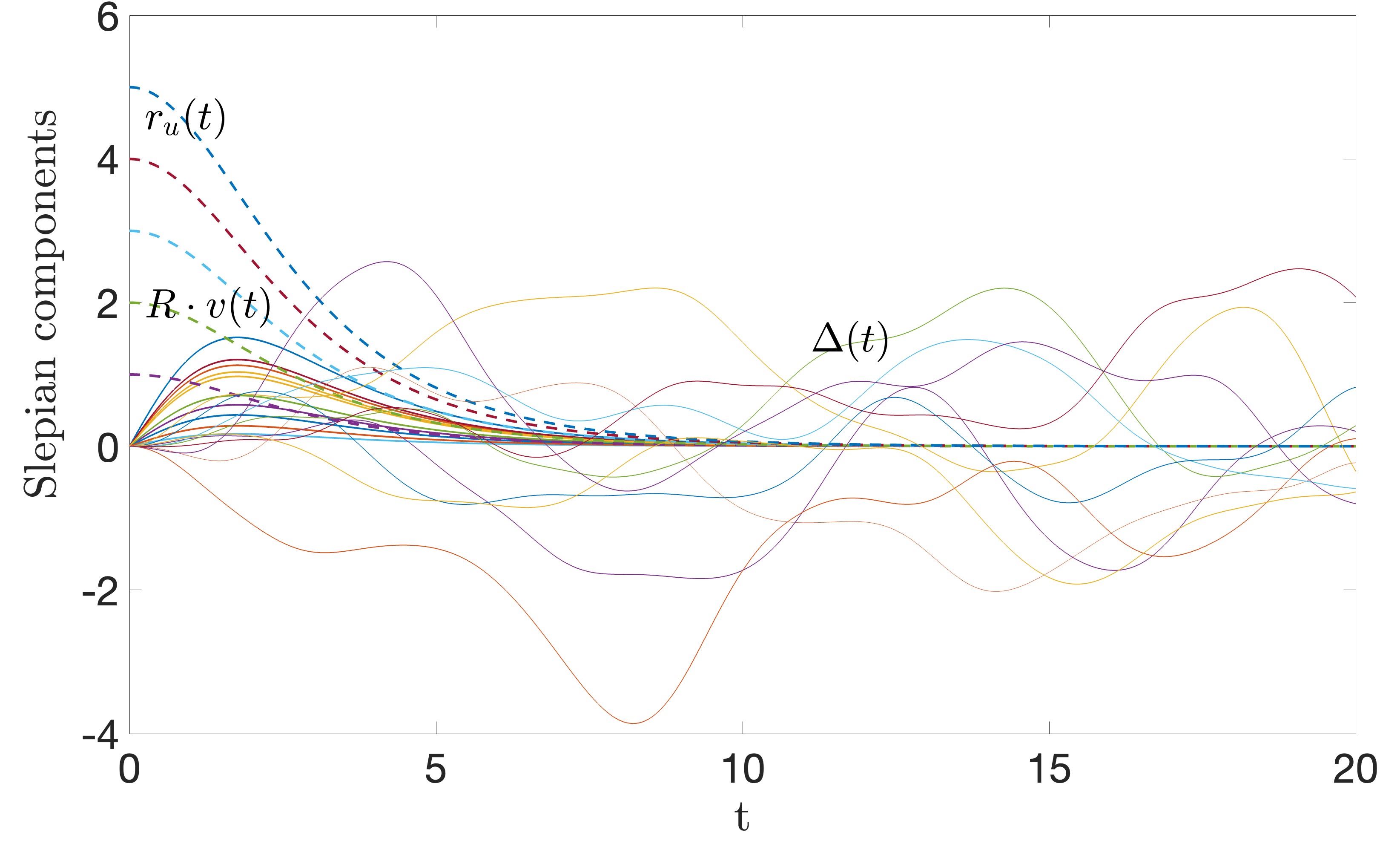}
\includegraphics[width=0.5\textwidth]{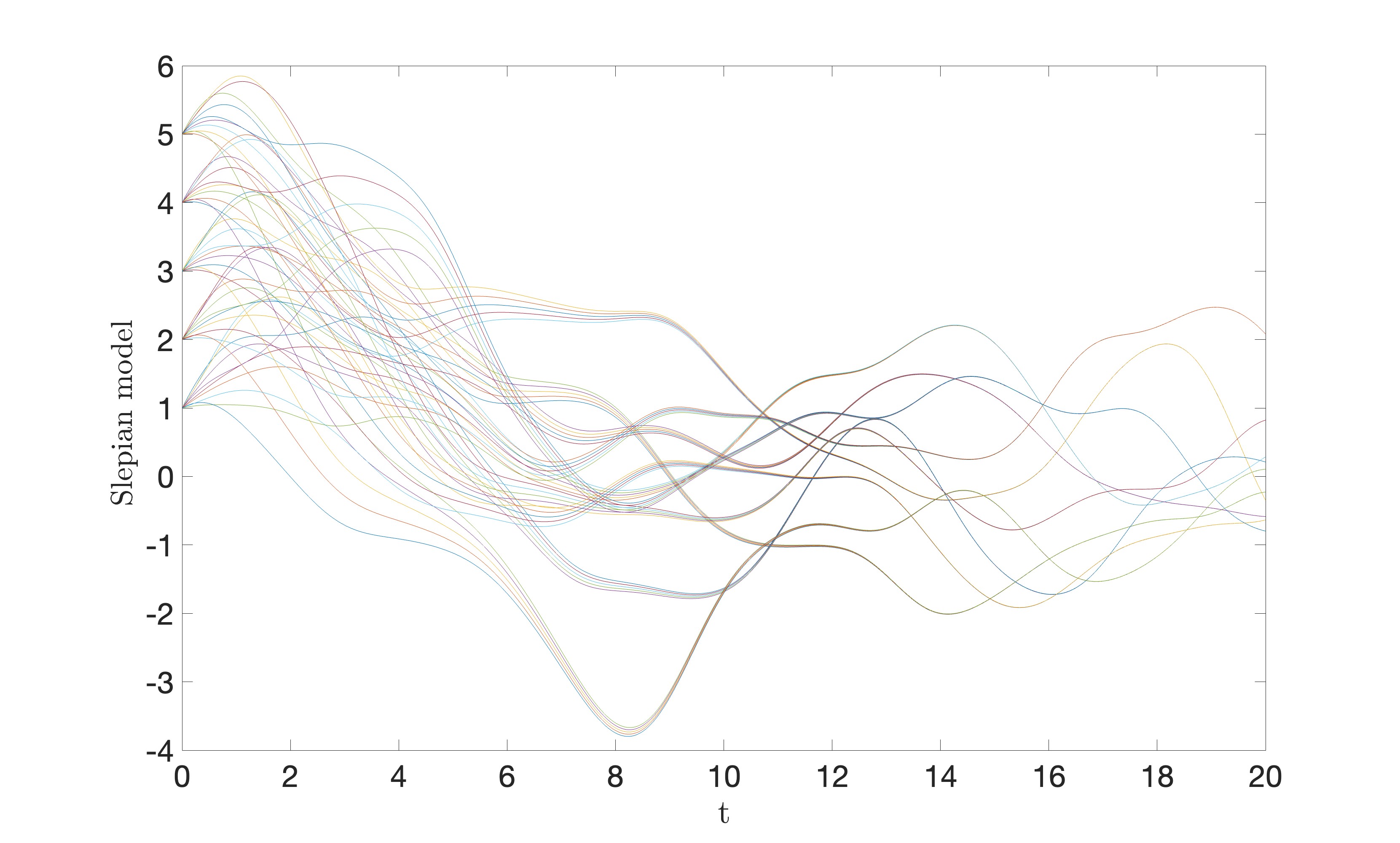}\\
\includegraphics[width=0.45\textwidth]{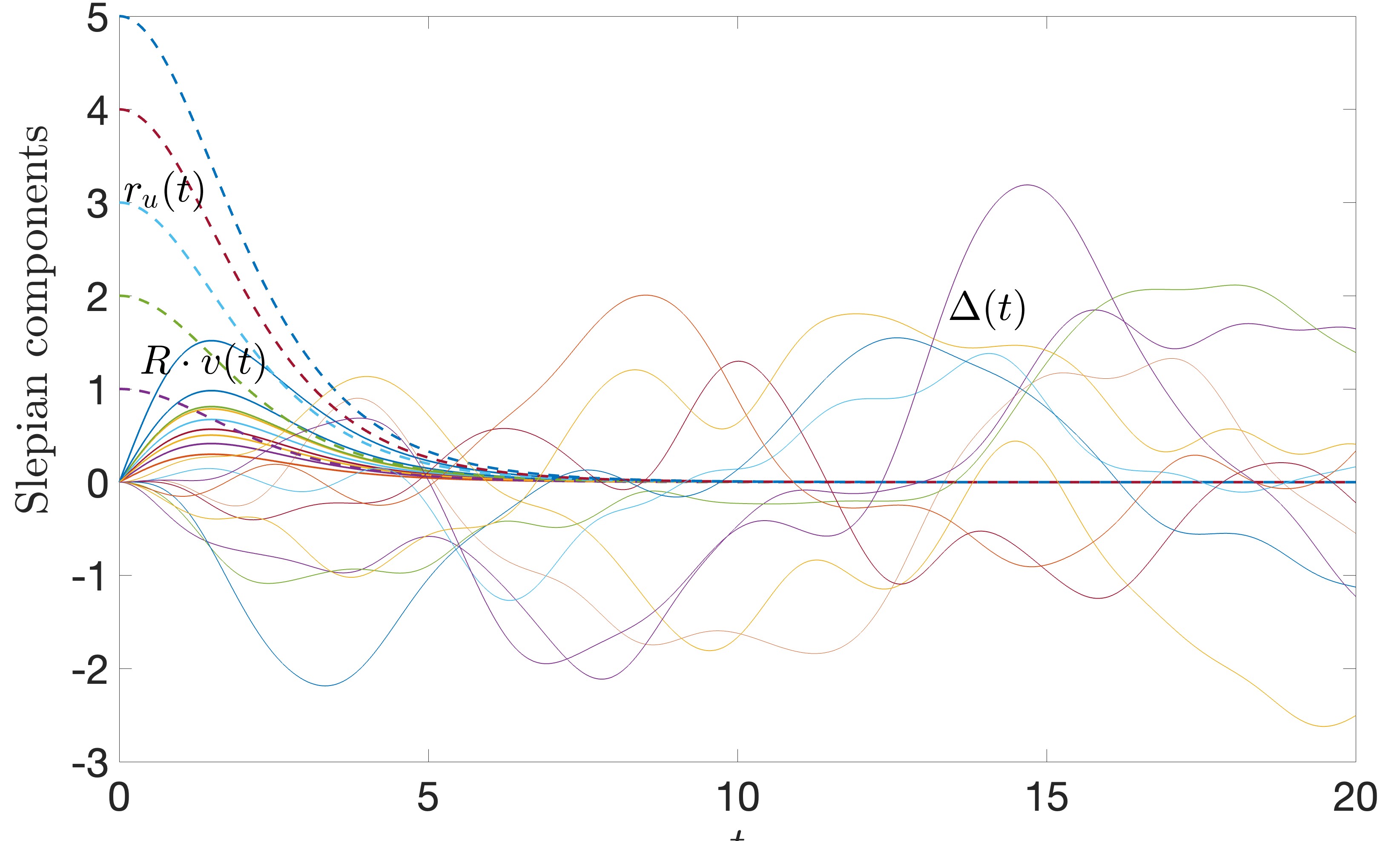}
\includegraphics[width=0.5\textwidth]{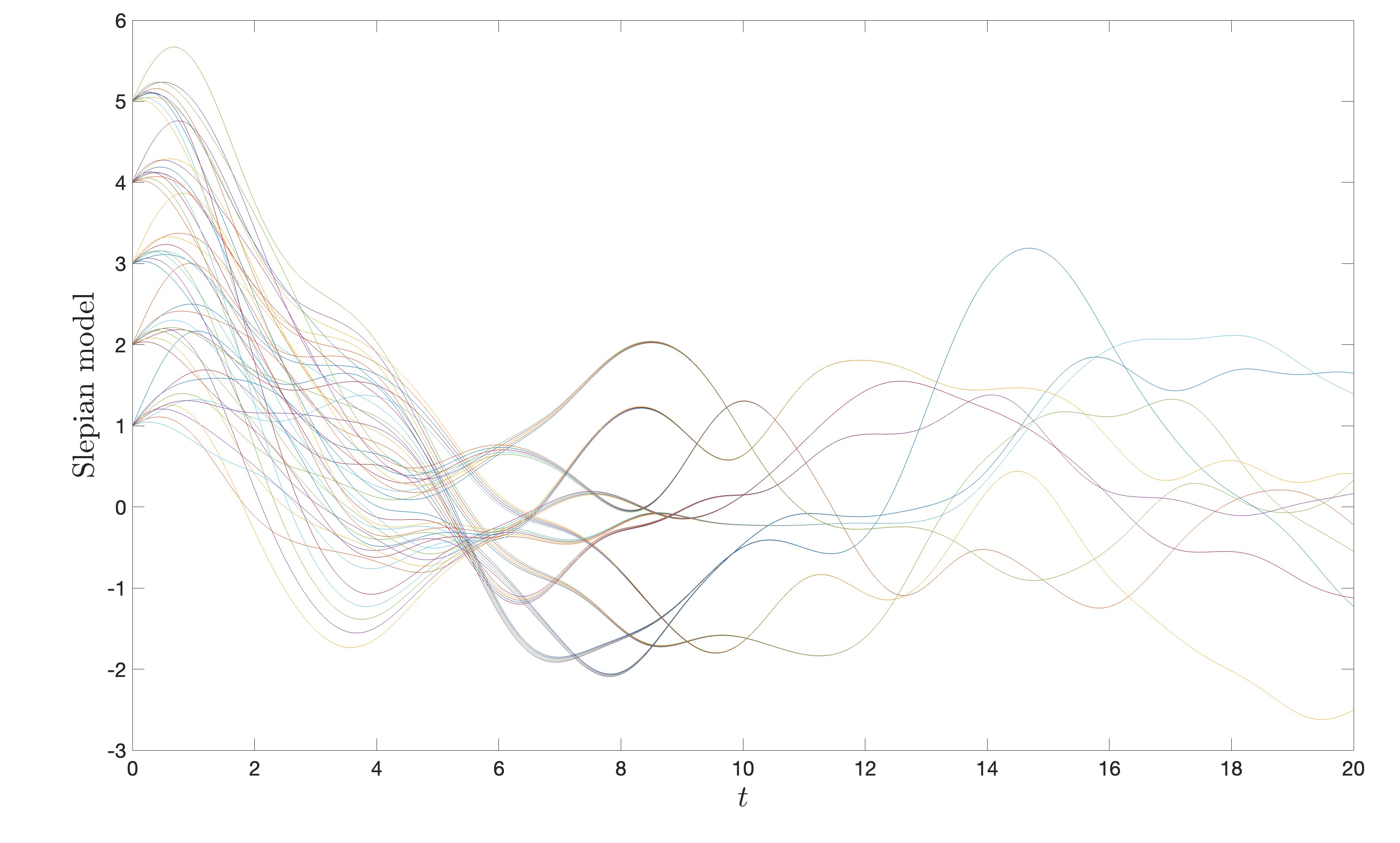}
\caption{\small The Slepian model for three Gaussian diffusions: 1-dimensional {\it (top)}, 2-dimensional {\it (middle)}, 3-dimensional {\it (bottom)}. {\it Right:} Three components of the decomposition; {\it Left:} Components put together to create ten samples from the Slepian model at five different levels.}
\label{fig:Slep}
\end{figure}

The Slepian process decomposes the behavior at the $u-$up crossing into three distinct parts. The first is fully deterministic and only depends on the level $u$: $r_u(t)=u \ r(t)$. The second is the deterministic function $-r'(t)/\sqrt{r''(0)}$, scaled by a standard Rayleigh random variable $R$. The last component is a non-stationary Gaussian stochastic residual process 
\begin{align*}
\Delta(t)=X(t)- X(0)\cdot r(t) - X'(0)\cdot r'(t)/r''(0),
\end{align*}
independent of $R$ and $u$.
In Figure~\ref{fig:Slep}, we present samples from these three components for Gaussian diffusion processes together with samples from the Slepian models combined from these components.

We also note that due to the symmetry of Gaussian processes, the Slepian model $\tilde X_u$ of the process $X$ at a $u$-down-crossing is given by
\begin{equation}
    \label{eq:downcr}
    \tilde X_u(t)\stackrel{d}{=}-X_{-u}(t), 
\end{equation}
where the relation denotes equality of the distributions of the two processes. Throughout this paper, unless otherwise stated, we will assume that all covariance functions are normalized such that $r(0)=1$. With the Slepian process introduced, we will move to the binary processes, which are used in the independent interval approximation method.

\section{Stationary and non-stationary switch processes} \label{sec:switch}
\noindent
The key idea of the Slepian-based IIA is to match the expected value of the clipped Slepian process to the expected value of a switch process. This section provides a short account of the necessary results for this matching and the results that later give the probabilistic motivation for why the Slepian-based IIA is valid. For an overview of the switch process, see \cite{JAPbengtsson} and \cite{BengtssonP}. The switch process is defined by interlacing intervals of random length such that after the end of each random length, the process switches either from one to minus one or from minus one to one. The process is defined on the positive half line starting at the origin, and the initial value is determined by a Bernoulli random variable, such that the process starts at one with some probability $p_0$ and at minus one with $1-p_0$. The Bernoulli random variable is independent of the intervals, which are denoted by $T_i^-$, and  $T_i^+$, for $i \in \mathbb N$ for the time spent in state minus one and one, respectively. Additionally, it's assumed that $T_i^+$ and $T_i^-$ sequences of independent identically distributed (iid) random variables which are also mutually independent.

\begin{figure}[t!] 
\begin{center}
    \includegraphics[width=0.75\textwidth]{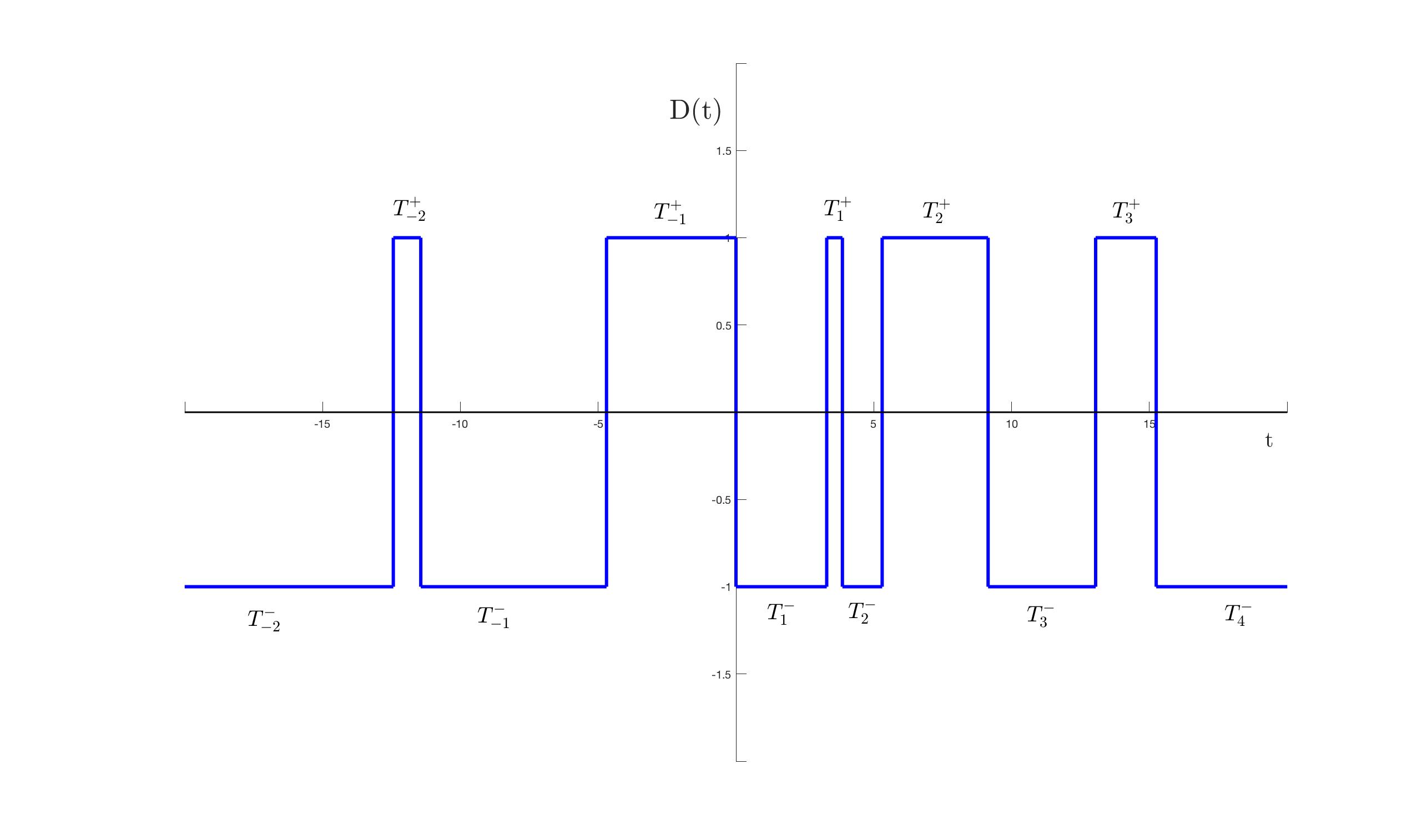} \vspace{-0.5cm}
\end{center}
\caption{\small The definition of the switch process attached to the origin, with the initial state $\delta=-1$.}
\label{SwitchLine}
\end{figure}

The resulting plus-minus process is denoted by $D(t)$, $t>0$. However, $D(t)$ can be extended to the negative real line by attaching another realization at zero that starts from $-\delta$ in the negative direction. In Figure~\ref{SwitchLine}, we illustrate a trajectory of such a process, which starts at minus one. 
Before presenting some elementary results needed for this paper, we need to establish the notation used for the Laplace transform since this is the main mathematical tool used. Let $\Psi(s)$ be the Laplace transform of a probability distribution defined by the cdf $F$ with support on the positive half line. The Laplace transforms of $T_i^+$, $i \in \mathbb N$ are denoted by $\Psi_+(s)$ and $\Psi_-(s)$.

In this work, we utilize Proposition 1 from  \cite{BengtssonP} on the Laplace transform of the expected value function of the switch process. 
\begin{proposition}
\label{OneDim}
Let $D(t)$ be a non-stationary switch process then the Laplace transform of $P_\delta(t)=\Prob{\left(D(t)=1 \big | \delta\right)}$, $t>0$, is given by
\begin{align*}
\mathcal L (P_\delta)(s)&=\frac{1-\Psi_{+}(s)}{s(1-\Psi_+(s)\Psi_-(s))} \cdot
\begin{cases}
1&;\delta=1,\\
\Psi_-(s)&; \delta=-1,
\end{cases}
\end{align*}
Additionally, the Laplace transform of the expected value function $E_\delta(t)=\Ex{( D(t)\big |\delta)}$, $t>0$ is given by 
\begin{align*}
\mathcal L (E_\delta)(s)&= \frac{\Psi_-(s)-\Psi_+(s)+\delta (1-\Psi_-(s))(1-\Psi_+(s))}{s(1-\Psi_+(s)\Psi_-(s))}.
\end{align*}
\end{proposition}
Focusing on the expected value functions to characterize this process is natural since $D(t)$ is not stationary. The reason for this lack of stationery is the attachment of the switch at zero. 

In the ordinary IIA, a stationary switch process is used to approximate the clipped process, and matching is done through the covariance function. A stationary switch process can be constructed from the non-stationary version by delaying the process forward and backward around zero. Hence, the origin will be placed in some interval $[-B, A]$, and two independent non-stationary switch processes are attached at each end of this interval. The value of the process on this interval is chosen by a binary random variable $\delta$, that takes the value one with probability $\mu_+/(\mu_++\mu_-)$ and minus one with probability $\mu_-/(\mu_++\mu_-)$. The main question is now how to find the distributions of $A$ and $B$ such that the delay results in a stationary process. The answer follows from well-known results in renewal theory see, for example, \cite{Cox1962}, and the densities, conditional on $\delta$, are given by 
\begin{align*}
    f_{A\vert \delta}(t)=f_{B\vert \delta}(t)= \frac{1-F_\delta(t)}{\mu_\delta}.
\end{align*}
Where $F_\delta$ is the cdf of $T_+$ or $T_-$ depending on if $\delta$ takes the value one or minus one, respectively. 
Throughout the rest of the paper, we will denote this stationary switch process by $\tilde{D}(t)$, $t \in \mathbb{R}$. 
The characteristic associated with the stationary switch process is the covariance function. We therefore present a shortened version of Proposition 2 from \cite{BengtssonP}, which gives us the covariance function of $\tilde{D}(t)$ in the Laplace domain. 
\begin{proposition}
\label{stat}
Let $\tilde{D}(t)$ be a stationary switch process constructed from the interval distributions $T_i^+$ and $T_i^-$ for $i\in\mathbb{N}$ with the expected values $\mu_+$ and $\mu_-$ respectively. Define $\tilde P_\delta(t)=P(\tilde D(t)=1|\tilde D(0)=\delta)$, then the process $\tilde{D}(t)$ is uniquely characterized by $\tilde{P_\delta}(t)$ which have the following form in the Laplace domain
\begin{align*}
\mathcal L (\tilde P_\delta)(s)&=\frac{1}{s}
\begin{cases} \displaystyle
1-\frac 1 {\mu_+ s} 
\frac{\left(1-\Psi_+(s)\right)\left(1-\Psi_-(s)\right)}{1-\Psi_+(s)\Psi_-(s)} 
&;\delta=1,\vspace{2mm}\\
\displaystyle
\frac 1 {\mu_- s} 
\frac{\left(1-\Psi_+(s)\right)\left(1-\Psi_-(s)\right)}{1-\Psi_+(s)\Psi_-(s)}&; \delta=-1.
\end{cases}
\end{align*}
The covariance of $\tilde{D}(t)$, $R(t)=\Cov{(\tilde D(u), \tilde D(t+u))}$ is given by 
\begin{align}
 \label{eq:CovExp}
R(t)=\frac{2}{\mu_++\mu_-}
 \left(
 \tilde P_{+}(t)\mu_+-\tilde P_{-}(t)\mu_-+
\mu_+ 
 \frac{\mu_--\mu_+}{\mu_++\mu_-}
 \right).
\end{align}
and its Laplace transform is given by 
\begin{align} \label{eq:LapR}
\mathcal L (R)(s)=\frac{4}{s\left(\mu_++\mu_-\right)}\left(\frac{\mu_+\mu_-}{\mu_++\mu_-}- \frac1s \frac{(1-\Psi_+(s))(1-\Psi_-(s))}{1-\Psi_-(s)\Psi_+(s)} \right).
\end{align}
\end{proposition}

We now have a collection of characteristics in the Laplace domain, $\Psi_+(s)$, $\Psi_-(s)$, $\mathcal{L}(E_+)(s)$, $\mathcal{L}(E_-)(s)$ and $\mathcal{L}(R)(s)$. The main goal now is to retrieve $\Psi_+(s)$ and $\Psi_-(s)$ from either the covariance or the expected value functions. It was shown in \cite{BengtssonP} that $\Psi_+(s)$ and $\Psi_-(s)$ are obtainable from the derivatives of the expected value function, and this results in the following expression
\begin{align} \label{eq:PsiExp}
    \Psi_+(s)=\frac{\mathcal L (E_+')(s)}{\mathcal L (E_-')(s)-2}, \ \ 
    \Psi_-(s)=\frac{\mathcal L (E_-')(s)}{\mathcal L (E_+')(s)+2}.
\end{align}
However, if one wants to use the above relations to match the expected value functions of the clipped Slepian process with the crossing distributions,  one needs to show that the above functions are completely monotone so that they indeed correspond to the Laplace transform of probability distributions. This is not easy to demonstrate directly from the definition. However, Theorem 1 in \cite{BengtssonP} gives sufficient conditions on the derivatives of $E_+$ and $E_+$, which can be used to solve this problem. 

\begin{theorem} \label{th:2}
Let $D(t)$ be a non-stationary switch process with the expected value functions $E_+(t)$ and $E_-(t)$. Then the following conditions are equivalent
\begin{itemize}
\item[\it i)] Functions $-E_+'(t)$ and $E_-'(t)$ are non-negative,\vspace{2mm}
\item[\it ii)] $T_+$ and $T_-$ have the stochastic representations
    \begin{align*}
        T_+\stackrel{d}{=} X+\sum_{k=1}^{\nu_\alpha-1} Y_k,  \ \ \ \
        T_-\stackrel{d}{=} Y+\sum_{k=1}^{\nu_\beta-1} X_k, 
    \end{align*}
    \end{itemize}
    where the sums are considered zero whenever summing over the empty set.
    Additionally, $X_i$'s and $Y_i$'s mutually independent and iid, with the densities 
    \begin{align*}
        f_X(t)=-\frac{1}{2\alpha}E_+'(t) \ \ \ \
        f_Y(t)=\frac{1}{2\beta}E_-'(t),
    \end{align*}
    respectively. $\nu_\alpha$ and $\nu_\beta$ are geometrically distributed with the parameters
    \begin{align*}
        \alpha=\frac{\mu_-}{\mu_++\mu_-} \ \ \ \ \beta=\frac{\mu_+}{\mu_++\mu_-},
    \end{align*}
    which also serves as the normalizing coefficients of $-E_+'$ and $E_-'$. 
\end{theorem}
\noindent
It should be mentioned that the monotonicity of the expected value functions is often easy to verify in practice. 
\vspace{-0.1cm}

\section{Clipping Slepian and Gaussian processes}
\label{sec:clipp}
\begin{figure}
\mbox{}\hspace{-1.1cm}\includegraphics[width=0.55\textwidth]{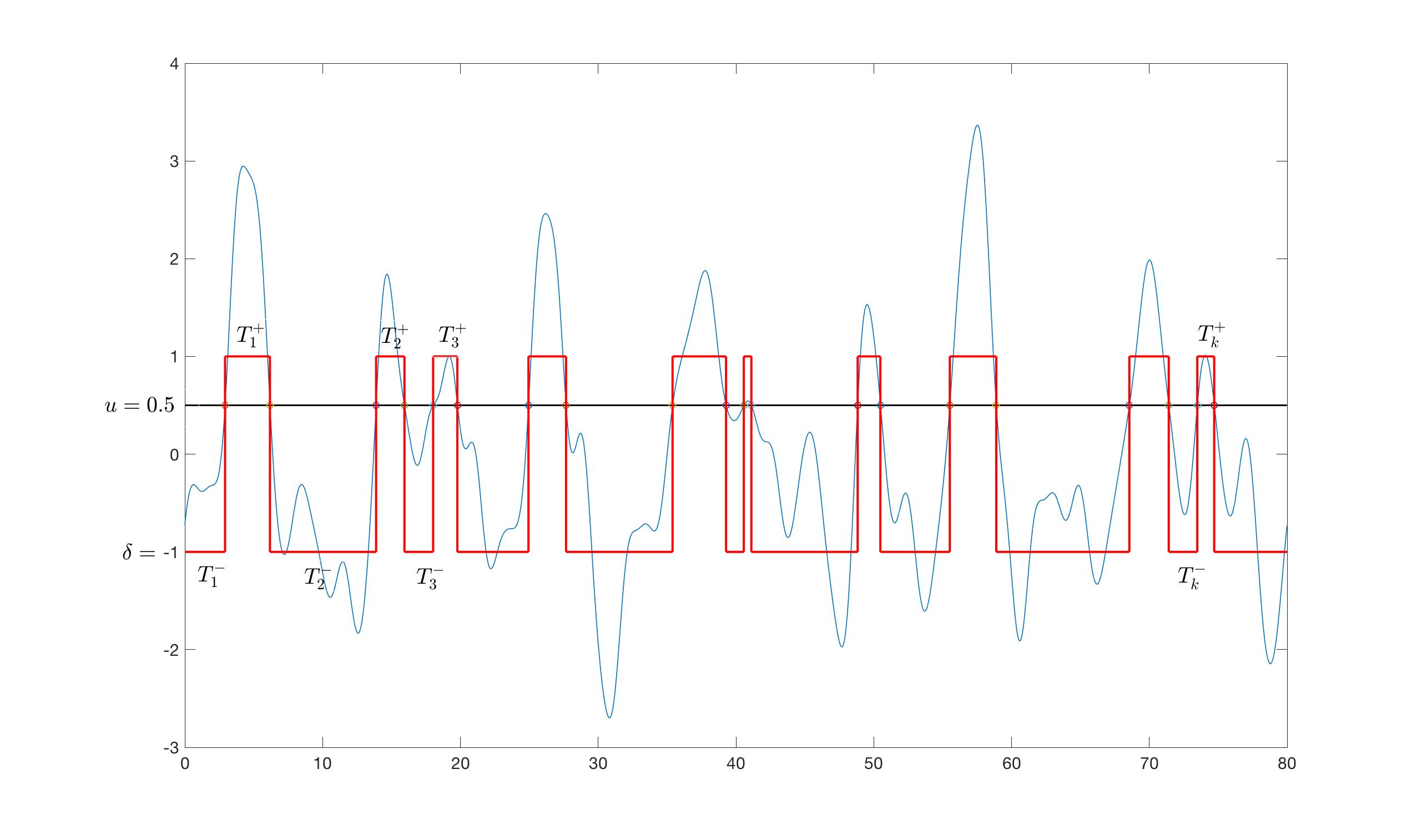}\hspace{-0.7cm}
\includegraphics[width=0.55\textwidth]{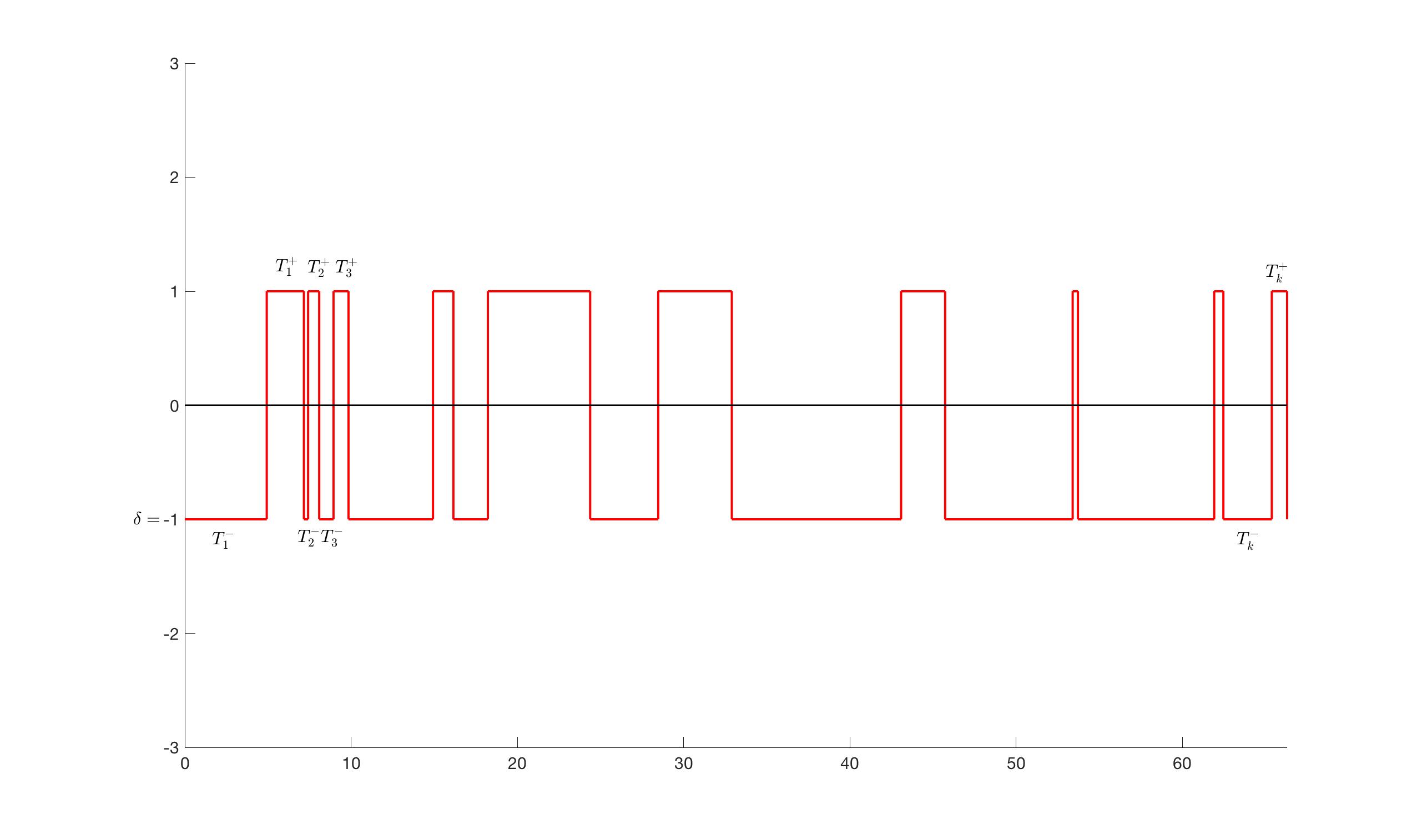}
\caption{\small The excursion intervals of a process $X(t)$ together with the corresponding clipped process (left) and a sample of the switch process under the IIA (right) with exponential distributions having the means matched those in the top graph.}
\label{swcl}
\end{figure}

\noindent
The length of $u$-level excursions of a smooth stationary process $X(t)$ forms a sequence of random variables. Additionally, if these excursions are marked with whether they fall below or above the level $u$, the process exhibits a structure similar to a switch process. However, the intervals for this process will not be independent, which is assumed for the switch process. We denote these interval lengths with $\bar{T}_i^+$ and $\bar{T}_i^-$, respectively, $i \in \mathbb Z$. More formally, we define the clipped process at level $u$ by 
\begin{align*}
    \tilde{D}_u(t)=sgn(X(t)-u).
\end{align*}
It is clear that $\tilde{D}_u(t)$ will inherit the stationarity property and contain information about the $u$-level excursions of $X(t)$ both above and below the level. In Figure \ref{swcl}, the idea of clipping a process is illustrated for a Gaussian process clipped at $u=0.5$. 
In the ordinary IIA approach, the covariance of the process $X(t)$, $r(t)$, is matched to the covariance of the clipped process. For the simplest case, the process is a zero mean Gaussian process clipped at $u=0$; we obtain the well-known and explicit arcsin formula
\begin{align*}
    R_0(t)=\frac 2 \pi \arcsin  r (t).
\end{align*}
For the non-zero level, greater care is needed. 
Let $F_X$ be the cdf of $X(0)$ and $\bar F_X=1-F_X$. 
Assuming continuity of the distribution of $X(t)$, the covariance of the clipped stationary process can be written as follows 
\begin{align*}
    R_u(t)&=\Prob{(X(t)>u | X(0)>u)}\bar F_X(u) +\Prob{(X(t) \le  u|X(0) \le u)}F_X(u)
    \\
&-\Prob{(X(t)\le u| X(0)>u)}\bar F_X(u) \Prob{(X(t)>u|X(0)\le  u)}F_X(u) - \left(1-2F_X(u)\right)^2
\\
&=2\left( \Prob{(X(t)>u| X(0)>u)}\bar F_X(u)+\Prob{(X(t)\le u|X(0)\le u)}F_X(u)+2F_X(u)\bar F_X(u)\right).
\end{align*}

If we assume that the process $X(t)$, $t\in \mathbb R$ is distributionally symmetric in value, i.e. it has the same distribution as $-X(t)$, $t\in \mathbb R$, then for $A_u(t)=\Prob{(X(t)\le u,X(0)\le u)}$ we have
\begin{align*}
R_u(t)=2\left(A_{-u}(t)+A_u(t)+2F_X(u)F_X(-u)\right).
\end{align*}
For a zero mean Gaussian process $X(t)$ with the covariance $r(t)$  and if we set $\tilde u=u/\sqrt{r(0)}$ and $\rho_t=r(t)/r(0)$: 
\begin{align*}
A_u(t)&=\Prob(\rho_tZ+\sqrt{1-\rho_t^2}Y<\tilde u, Z<\tilde u)\\
&=\frac{1}{\sqrt{2\pi}}\int_{-\infty}^{\tilde u}e^{-\frac{z^2}2} {\Phi}\left(\frac{\tilde u- \rho_tz}{\sqrt{1-\rho_t^2}}\right) ~dz=B(\tilde u, \rho_t),
\end{align*}
where $Z$ and $Y$ are independent standard normal variables, and the last equation serves as a definition of $B$. Consequently, we obtain the covariance of the clipped processes in terms of the covariance of the underlying Gaussian process
\begin{align}
\label{eq:autcov}
R_u(t)&=2\left(B\left(\frac{u}{r(0)},\frac{r(t)}{r(0)}\right)+B\left(\frac{-u}{r(0)},\frac{r(t)}{r(0)}\right)+\Phi\left(\frac{u}{r(0)}\right)\Phi\left(\frac{-u}{r(0)}\right)
\right).
\end{align}
Since the clipped process of a stationary process is stationary, it is natural to characterize the clipped process by its covariance function.  

However, in the IIA approach, one has to match two characteristics corresponding to the distributions of the excursion above and below the level $u$. For this purpose, the Slepian process attached at $u$-crossing (up or down) constitutes a much better alternative than the clipped stationary Gaussian process. The reason for this is that we do have two Slepian processes at a crossing of $u$. One is for the $u$-upcrossing and the other for the $u$-downcrossing. In our approach, we clip two Slepian processes given in \eqref{eq:Slepian} and \eqref{eq:downcr}
\begin{align}
\label{eq:clippedSlepianup}
D_u^+(t)&={\rm sgn}(X_u(t)-u)\\
\label{eq:clippedSlepiandown}
D_u^-(t)&={\rm sgn}(\tilde X_{u}(t)-u).
\end{align}
We note that by  \eqref{eq:downcr}, $D_u^-\stackrel{d}{=}-D_{-u}^+$.

In the next proposition, the expected value of the clipped process is presented. Later, this function will be used to obtain the approximated excursion distribution by matching it with a non-stationary switch process. Hence, this is an important equation for the Slepian-based IIA. 
\begin{proposition}
\label{th:expctSlep}
Let $D^+_u(t)$ be the clipped Slepian process at an $u$-level up crossing of a Gaussian process, $u\in \mathbb R$, with the covariance function $r(t)$ twice continuously differentiable.
Then $E_u^+(t)=\Ex(D_u^+(t))$ is given by
\begin{align*}
    E_u^+(t)&=1- 2 \Phi\left( u
    \frac{1-r(t)}{\sqrt{1-r^2(t)+\frac{r'^2(t)}{r''(0)}}}
    \right) 
    -2 \frac{1}{\sqrt{-r''(0)}}\frac{r'(t)}{\sqrt{1-r^2(t)}} \exp\left(-\frac{u^2}{2}\frac{1-r(t)}{1+r(t)} \right) 
    \\
    &\times \Phi\left(\frac{-u}{\sqrt{-r''(0)}} \sqrt{\frac{1-r(t)}{1+r(t)}}\frac{r'(t)}{\sqrt{1-r^2(t)+\frac{r'^2(t)}{r''(0)}}}\right).
\end{align*}
Let  $E_u^-(t)=\Ex(D_u^-(t))$, then  
\begin{align*}
    E_u^-(t)=-E_{-u}^+(t). 
\end{align*}
Moreover, 
\begin{align*}
    \lim_{t\rightarrow 0^+}E_u^+(t)&=1,\ \ \ \, \lim_{t\rightarrow 0^+}E_u^-(t)=-1,\\
    \lim_{t\rightarrow \infty }E_u^+(t)&=\lim_{t\rightarrow \infty}E_u^-(t)=1-2\Phi\left(u \right).
\end{align*}
\end{proposition}
The proof requires some technical lemmas, and these, together with the proof, have been relegated to the Appendix. While the expression of the clipped Slepian process might appear to be intimidating, it simplifies significantly for the zero level. For $u=0$, it reduces to
\begin{align*}
     E_0^+(t)=- \frac{1}{\sqrt{-r''(0)}} \frac{ r'(t)}{\sqrt{1-r(t)^2}}. 
\end{align*}
In the next section, the characteristics of this section, such as the covariance of the clipped process and expected value functions of the clipped Slepian process, will be used to obtain approximations of the excursion distributions.

\section{IIA for a non-zero crossing level}
\label{sec:IIA}
\noindent
Intuitively, the IIA simply means that the characteristics of a clipped process are matched to the characteristics of a switch process, either the stationary one as in the ordinary IIA case or the non-stationary one in the Slepian-based IIA. Both approaches imply that the dependency between the $\bar{T}_i^+$'s and $\bar{T}_i^-$'s of the clipped process is presumed to be negligible so that independence can, for practical purposes, be assumed. After matching the characteristics, the relation between the characteristics and the switching times of the switch processes is utilized to obtain the switching time distributions, which serve as an approximation of the excursion time distributions. We start with an overview and treatment of the ordinary IIA and then introduce the Slepian-based IIA, which extends to non-zero levels of the work done in \cite{bengtsson2024slepian}.

\subsection{The stationary IIA}
\noindent
Initially, when zero-level crossings were considered, the stationary versions of the processes were matched through their covariance functions, see \cite{BrayMS}. 
It was shown in \cite{BengtssonP} that the Slepian-based IIA coincides with the ordinary IIA for zero-level crossings, i.e. that the two approaches are equivalent. 
However, the non-zero crossing case is more complicated and has not been extensively studied. 
Using characteristics of a clipped process and matching it to characteristics of a stationary switch process was introduced in \cite{Sire2007} and further investigated in \cite{Sire2008}. Here, we present a more rigorous introduction of the IIA at the non-zero crossings for the stationary case, which is needed to identify the gaps in the methodology that will be later addressed through the Slepian-based approach.

Following the previous sections, we have two stationary processes, the switch stationary one $\tilde D(t)$, $t\in\mathbb R$, and the clipped stationary Gaussian one $\tilde D_u(t)={\rm sgn}(X(t)-u)$, $t\in\mathbb R$. It is clear that the stationary switch process $\tilde D(t)$ is simpler in structure and can be utilized for approximate analysis of $\tilde D_u(t)$ if the two are matched in one or the other way. In the more straightforward, zero-level/symmetric case, the two processes can be matched by their covariance structure, i.e., we call $\tilde D(t)$ the independent intervals approximation of $\tilde D_0(t)$ if the covariances of the two are matched, i.e. 
\begin{align*}
R(t)=R_0(t).
\end{align*}
For the case of $u=0$, it simply means that the  (identical) distribution of $T_i^+$ and $T_i^-$ is given by the cdf $F$, the following relation have to be satisfied
\begin{align*}
\mathcal L\left(\frac 2 \pi \arcsin \rho_t\right)(s)
=
\frac  1 s \left(1-\frac{2}{s\mu} \frac{1-\Psi(s)}{1+\Psi(s)}\right).
\end{align*}
Through this relation, the distribution given by $\Psi$ matches the covariance of the original process. Conversely, if the distribution of the excursion times is specified, then, through the above relation, we can obtain the covariance of the original process. This one-to-one correspondence (if such one exists) formally defines the IIA for zero-level crossings. 
Of course, it is not obvious that the relation leads to a valid distribution, which has been discussed in \cite{bengtsson2024slepian}.
In this case, one needs also match also $\mu$, which can be done by the reciprocal of the crossing intensity given by the Rice formula
$$
\lambda_0=\frac{\Ex\left(|X'(0)|\right)}{\sqrt{2\pi r(0)}}=\frac{1}{\pi}\sqrt{\frac{-r''(0)}{r(0)}},
$$
so that the matching equality is 
 $$
 \mu=\frac{1}{\lambda_0}=\pi \sqrt{\frac{r(0)}{-r''(0)}},
 $$

We note that for the asymmetric case, i.e., when $u\ne0$, there are two distributions to be matched, $F_+$ for the excursions above the $u$-level and $F_-$ corresponding to these below the $u$-level. 
Thus, if we consider matching the covariance of the switch process and the clipped process, one needs one more relation to solve for both distributions. For that, different strategies could be taken.

Let us consider the matching through the covariance functions. First, the covariance $R_u$ of the clipped process is given in  \eqref{eq:autcov}. The relation between the covariance function of the stationary switch process and the interval distributions is given by Proposition \ref{stat}. By computing the derivative of  (\ref{eq:CovExp}) and using  and (\ref{eq:PsiExp}) we obtain the following 
\begin{align*}
\mathcal L (R'')(s)&=\frac{4}{\mu_++\mu_-}\cdot \frac{\Psi_+(s)+\Psi_-(s)-2\Psi_+(s)\Psi_-(s)}{1-\Psi_+(s)\Psi_-(s)}.
\end{align*}
As it can be seen from this relation, the roles of $\Psi_+$ and $\Psi_-$ cannot be distinguished solely based on the covariance function $R$ of the stationary switch process and thus from the covariance $R_u$ of the clipped process that matches $R$.  
There is a need for an additional relation.

Before discussing this, let us point out how $\mu_+$ and $\mu_-$ can be matched. 
First, we have $\mu_++\mu_-$ can be matched by the reciprocal of the half of  the $u$-level crossing intensity $\lambda_u$ which is given by the Rice formula
$$
\lambda_u=\frac{
\Ex\left(|X'(0)|\right)}{\sqrt{2\pi r(0)}}\,
e^{-{u^2}/{(2r(0))}}
=\frac{1}{\pi}\sqrt{\frac{-r''(0)}{r(0)}}\,
e^{-{u^2}/{(2r(0))}}
$$
so 
\begin{equation}
\label{eq:sum}
\mu_++\mu_-=2\pi \sqrt{\frac{r(0)}{-r''(0)}}e^{{u^2}/{(2r(0))}}.
\end{equation}
It is easy to argue that the ratio of averages of spending above and below the level $u$ must be equal to the ratio of the probabilities of $X(0)$ being above and below $u$, yielding the second matching equation
\begin{equation}
\label{eq:ratio}
\frac{\mu_+}{\mu_-}=\frac{\Phi\left(-u^2/(2r'(0))\right)}{1-\Phi\left(-u^2/(2r'(0))\right)},
\end{equation}
where $\Phi$ is the cdf of a standard normal distribution.

In \cite{Sire2007}, where non-zero level crossing IIA has been attempted for the first time, and two matching characteristics have been proposed
\begin{align}
\label{eq:Au}
R_{<u}(t)&=\Ex ((u-X(t))^+(u-X(0))^+),\\
\label{eq:Nu}
N_{<u}(t) & = \Ex \left(N\left(\{0<s\le t; X(s)=u\}\right)\,\mathbb I_{X(0)<u}\right),
\end{align}
where $x^+$ is a positive part of a real number $x$, $\mathbb I_A$ is an indicator function of a set $A$, and $N(\cdot)$ is counting measure of a set. 

These two characteristics of the clipped process parallel equivalent characteristics of the switch process.
Let us introduce the switch counting process $\tilde N(t)$, $t>0$ that counts switches in $(0,t]$ for the stationary switch process $\tilde D$. Clearly, 
\begin{equation*}
\label{rep}
\tilde D(t)=(-1)^{\tilde N(t)+(1-\delta)/2},\, t\ge 0.
\end{equation*}
The quantities corresponding to \eqref{eq:Au} and \eqref{eq:Nu} are, respectively,
\begin{align*}
A_<(t)&=\Ex \left(\frac{1-\tilde D(t)}2\frac{1-\tilde D(0)}2 \right ),\\
 N_<(t) & = \Ex \left(\tilde N(t)|\delta=-1 \right).
\end{align*}
We have 
\begin{align*}
A_<(t)&=\frac{R(t)}4 +\frac{\mu_-^2}{\left(\mu_++\mu_-\right)^2}.
\end{align*}
Thus, by utilizing Theorem~\ref{stat} and  Proposition~\ref{prop:exp} given in the Appendix, we obtain the following Laplace transform relations
\begin{align}
\label{eq:ALap}
\mathcal{L}\left(A_<\right)(s)&=\frac{1}{s\left(\mu_-+\mu_+\right)}
\left(\mu_- - \frac1s \frac{(1-\Psi_+(s))(1-\Psi_-(s))}{1-\Psi_-(s)\Psi_+(s)} \right),
\\
\label{eq:N<Lap}
\mathcal L\left(N_{<}\right)(s)&=\frac{
 1
 }{s^2\mu_-}
  \frac{ \left( 1 +
  \Psi_{+}(s)\right) \left(1-\Psi_{-}(s)\right)}{1-\Psi_+(s)\Psi_-(s)}.
\end{align}
Given the characteristics $A_<$ and $N_<$, these two equations allow the identification of both $\Psi_+$ and $\Psi_-$. 
\begin{remark}
We note that the formula for the Laplace transform of $N_<$ coincides with that in \cite{Sire2007} and \cite{Sire2008}, and the transform of $A$ is also matching presented there if one accounts matching for $\mu_+$ and $\mu_-$ given in \eqref{eq:sum} and \eqref{eq:ratio}.
\end{remark}

In the stationary IIA approach, which is presented above, one challenge is to evaluate explicitly $N_{<u}$ for the clipped process. Moreover, it is not apparent that matching through \eqref{eq:ALap} and \eqref{eq:N<Lap} leads to valid distribution functions given through $\Psi_+$ and $\Psi_-$.

\subsection{The Slepian-based IIA}
\noindent
The core idea of the Slepian-based IIA is to use the Slepian process and Theorem~\ref{th:expctSlep} to obtain approximated excursion distributions as the distributions  $\hat{T}_+$ and $\hat{T}_-$ obtained through (\ref{eq:PsiExp}) as follows
\begin{align}
\label{eq:SlIIA}
\Psi_{\hat{T}_+}(s)=\frac{\mathcal L {(E_u^+}')(s)}{\mathcal L  {(E_u^-}')(s)-2}   \ \ \ \ 
\Psi_{\hat{T}_-}(s)=\frac{\mathcal L {(E_u^-}')(s) }{\mathcal L {(E_u^+}')(s)+2}.
\end{align}

For the equations to correspond to valid probability distributions, they need to be completely monotone. To verify it for a given pair of ${E_u^+}'$ and ${E_u^-}'$ is difficult at best. 
This is where the utility of Theorem \ref{th:2} becomes apparent.
Since it provides easy-to-check conditions on ${E_u^+}'$ and ${E_u^-}'$ such that they correspond to valid probability distributions, the next theorem formulates this in the context of the clipped Slepian process.
\begin{theorem}
    \label{th:part_char}
    Let $X$ be a stationary Gaussian process with the twice continuous differentiable covariance function $r(t)$ $t\in \mathbb{R}$. If for a given level $u$ the functions $E_u^+$ and $E_u^-$ from Proposition~\ref{th:expctSlep} are monotone, then the approximated excursion distributions above and below the level $u$ in the Slepian based IIA has the stochastic representation
     \begin{align*}
        \hat{T}_+\stackrel{d}{=} X+\sum_{k=1}^{\nu_\alpha-1} Y_k, \ \ \ \ 
        \hat{T}_-\stackrel{d}{=} Y+\sum_{k=1}^{\nu_\beta-1} X_k, 
    \end{align*}
    where sums over the empty set are zero when the sum is empty. $\nu_\alpha$ and $\nu_\beta$ are geometrically distributed with parameters
\begin{align*}
    \alpha=\Phi\left({u}\right) \ \ \ \ \beta=1-\Phi\left({u}\right)
    \end{align*}
    respectively, and the variables $X_i$'s and $Y_i$'s have densities
    \begin{align*}
        f_X(t)=-\frac{1}{2\alpha}{E_u^+}'(t) \ \ \ \
        f_Y(t)=\frac{1}{2\beta}{E_u^-}'(t).
    \end{align*}
    The sequences $X_i$ and $Y_i$ are iid and independent of $\nu_\alpha$ and $\nu_\beta$ which are also mutually independent. 
\end{theorem}
As previously mentioned, for the zero level case, it was shown in \cite{bengtsson2024slepian} that the ordinary IIA and the Slepian-based IIA coincide. This does not seem to be true for the non-zero level crossings since $A_{<u}$ and $N_{<u}$ of the clipped process do not seem related to $E_u^+$ and $E_u^-$ in the same way as $A_<$ and $N_<$ are related to $E_+$ and $E_-$ for the switch process. We advocate for the Slepian-based IIA given through  \eqref{eq:SlIIA} since $E_u^+$ and $E_u^-$ explicitly given in Theorem~\ref{th:expctSlep}, while $A_{<u}$, that can be obtained from  \eqref{eq:autcov}, involves integration formulas and obtaining a formula for $N_{<u}$ is a nontrivial open problem, yet to be solved.

\section{Non-zero level persistency coefficients for Gaussian diffusion}
\label{sec:diff}
\noindent
This section presents an example of the IIA for non-zero level crossing for the Gaussian diffusion process. We focus on the $2-$dimensional diffusion process for two reasons. The first reason is that it is one of the few processes for which the persistency coefficient is analytically known for the zero level, see \cite{PoplavskyiSchehr2018}. The non-zero excursion has also been studied by, among others \cite{Sire2008} using the ordinary IIA. The second reason is its wide use in statistical physics, where the diffusion covariances were used in \cite{MajumdarSBC}, \cite{BrayMS}, and \cite{WongMairWalsworthCory2001}.

In particular, \cite{WongMairWalsworthCory2001} studied the case where a diffusing field starts from a random initial configuration.  This model is theoretically attractive due to the linear form of the diffusion equation in any dimension and experimentally because of its natural physical interpretation in various non-equilibrium systems. For these reasons, the model constitutes a mathematically convenient and physically important benchmark for any methodology aiming at a persistency approximation and assessment.

For the sake of this presentation, it is sufficient to know that the covariance of the $d$-dimensional diffusion is given by $\mbox{\rm sech}^{d/2}(t/2),~d\in \mathbb N $ and that for $d=2$ we have an explicit spectrum $\mbox{\rm sech}(\pi \omega)$. While simulating trajectories might seem adequate, it becomes impractical for higher levels of $u$ where these crossing events occur less frequently, necessitating the simulation of very long trajectories, which is computationally costly.

The first step in approximating the persistency coefficients in the Slepian-based IIA framework is to verify that the monotonicity conditions in Theorem~\ref{th:part_char} are satisfied. This is verified by inserting the covariance function and its derivatives of the two-dimensional diffusion processes into the equations of Theorem~\ref{th:expctSlep}. Following this, the expected value functions can be turned into cumulative distribution functions (CDFs).

\begin{figure}[!t]
\centering
    \begin{minipage}{0.5\textwidth}
        \includegraphics[width=0.9\textwidth]{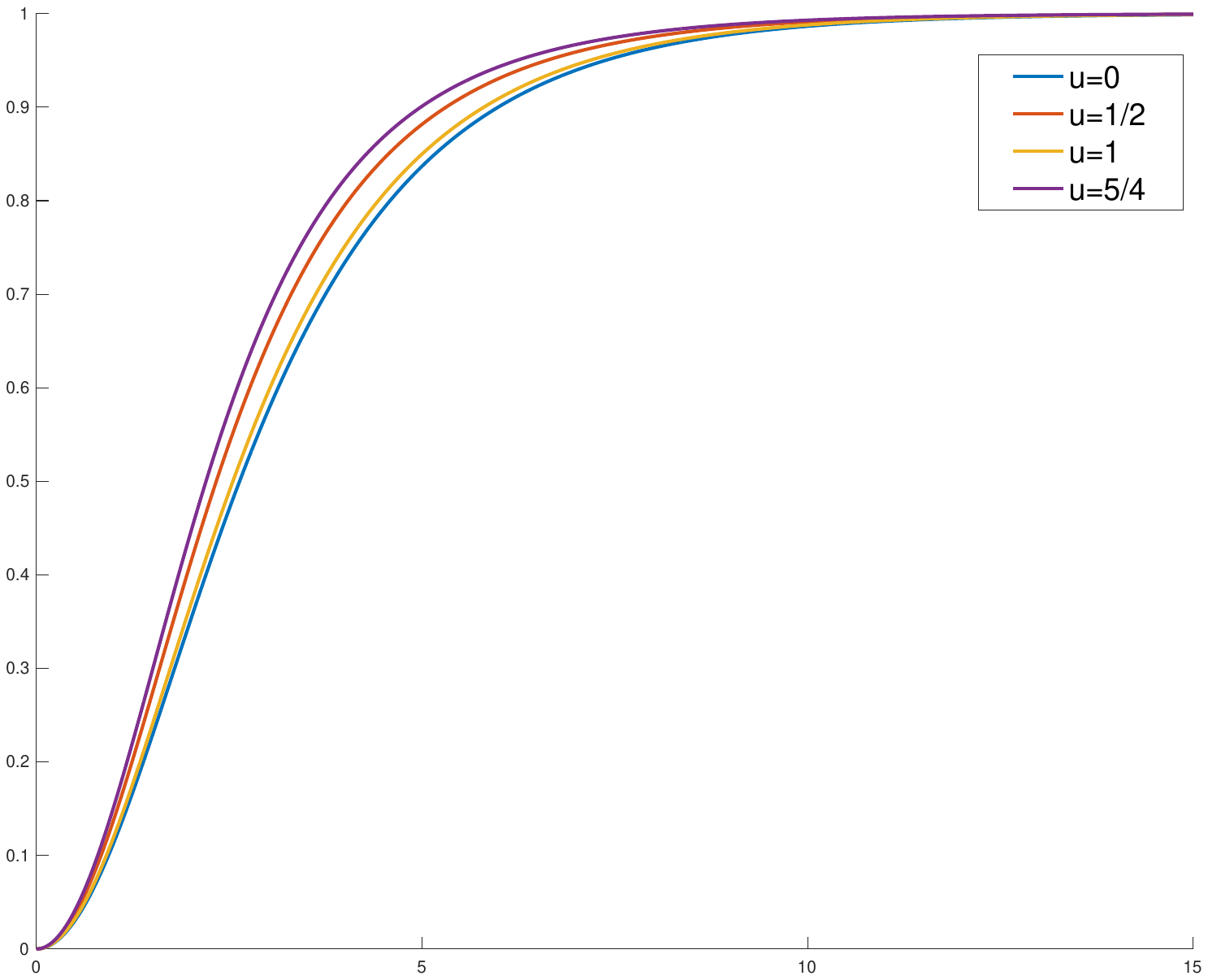}
    \end{minipage}\hfill
        \begin{minipage}{0.5\textwidth}
        \includegraphics[width=0.9\textwidth]{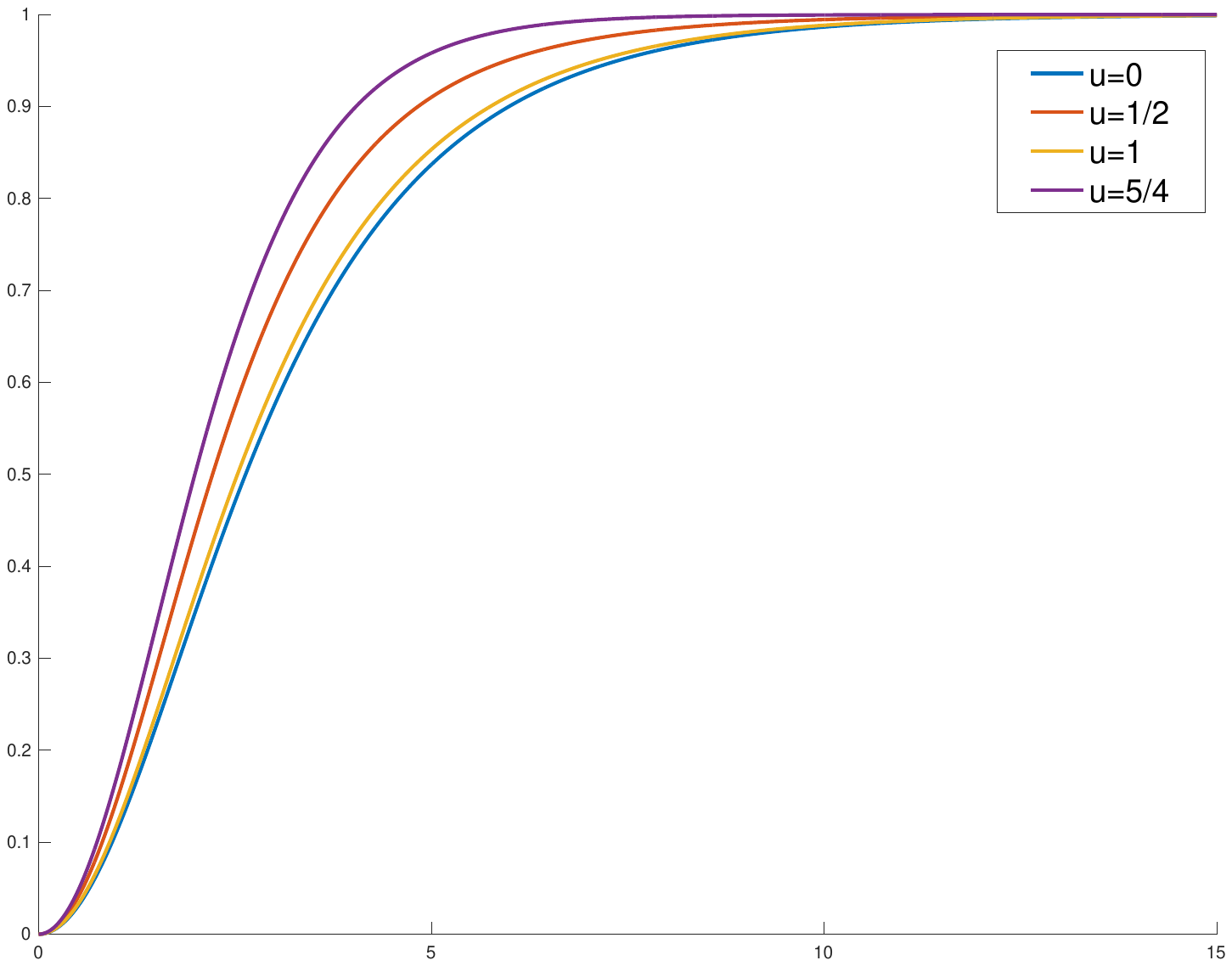}
    \end{minipage}\hfill
\caption{The expected values of the Slepian model based on the diffusion covariance, for an $u$ up-crossing {\it (left)} and an $u$ down-crossing {\it (right)}, for $u=0,1/2,1,5/4$.}
\label{fig:expgraph}
\end{figure}

Figure~\ref{fig:expgraph} illustrates the normalized expected value functions of the Slepian process for various levels of $u$. They correspond to the CDFs of the geometric divisors in the stochastic representation outlined in Theorem~\ref{th:2} for multiple levels of $u$.
It is evident from the graphs that the conditions of Theorem~\ref{th:part_char} are satisfied. However, it is important to note that $E_-$ loses its monotonicity properties for $u>5/4$. Hence, we focus on the crossing levels for which the monotonicity conditions in the previous theorem are met. To estimate persistence, we utilize samples from the approximated excursion distribution. These samples are obtained by numerical inversion sampling from the CDFs of the normalized expected value functions depicted in Figure~\ref{fig:expgraph}. The samples are then generated from the approximate excursion time distribution using the stochastic representation in Theorem~\ref{th:part_char}. However, the problem of estimating the persistency coefficients, $\theta_+$ and $\theta_-$, from these generated samples remains.

Suppose that the tails of $P(T_+\geqslant t)$ and $P(T_-\geqslant t)$ are of an exponential form, i.e., $P(T_+\geqslant t) \approx C \cdot e^{-\theta_+ t}$, for some large $t$, then the natural logarithm of $P(T_+\geqslant t)$ will be on the linear form $ln(C) - \theta_+ t$.  Hence, the persistency coefficients can be approximated via ordinary least squares (OLS) on the empirical survival function. While the notion of utilizing the empirical survival function for tail estimation is not new, having been introduced by \cite{KratzResnick} and \cite{SchultzeSteinebach}, it was employed in \cite{bengtsson2024slepian} for estimating persistence coefficients from generated samples. We obtain $10^6$ samples from the approximated excursion time distribution to estimate the persistence coefficients. Only the values below one-half of the survival function are utilized for the OLS approximation of $\theta_\pm$, as only the tail is of interest. This approach is repeated to obtain $10$ estimates, and approximate $95\%$ confidence intervals are computed and presented in Table~\ref{tab:Slepian}.

\begin{table}[h]
\centering
\caption{Slepian-based IIA persistency coefficient estimate}
\label{tab:Slepian}
\begin{tabular}{ccccc}
\toprule[0.9pt]
Crossing level, $u$ & Estimate, $\theta_+$ & Confidence interval & Estimate, $\theta_-$ & Confidence interval 
\\ 
\midrule[0.9pt]
0 & $0.1862$ & $\pm 0.000217$  & $0.1861$ & $\pm 0.000195 $  \\ 
$1/2$ & $0.2893$ & $\pm 0.000347$ & $0.1105$ & $\pm 0.000131 $  \\
$1$ & $0.4225$ & $\pm 0.000462 $ & $0.0591$ & $\pm 0.000061 $ \\
$5/4$ & $0.5001$ & $\pm 0.000574$ & $0.0411$ & $\pm 0.000044 $\\ 
\bottomrule[0.9pt]
\end{tabular}
\end{table}
For comparison, trajectories of the $2-$dimensional  Gaussian diffusion process were simulated using the WAFO-toolbox \citep{gitwafo}. A total of $10^5$ trajectories were simulated, each with a length of $10^6$. The persistence coefficient was estimated in the same manner as before, using OLS on the extracted crossing intervals. The estimates are presented in Table~\ref{tab:Wafo} along with approximate $95\%$ confidence intervals.

\begin{figure}[!t]
\centering
    \begin{minipage}{0.5\textwidth}
        \includegraphics[width=0.9\textwidth]{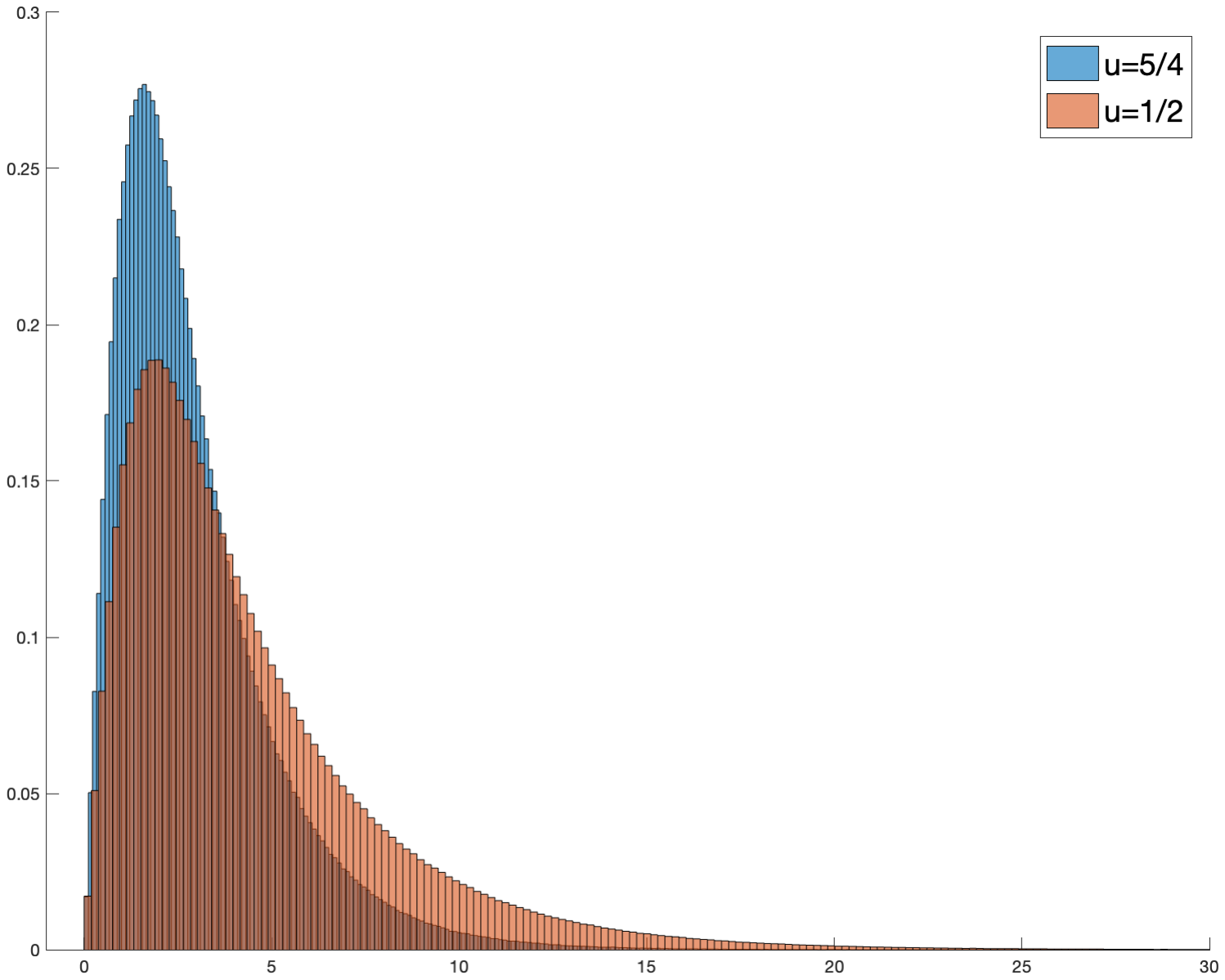}
    \end{minipage}\hfill
        \begin{minipage}{0.5\textwidth}
        \includegraphics[width=0.9\textwidth]{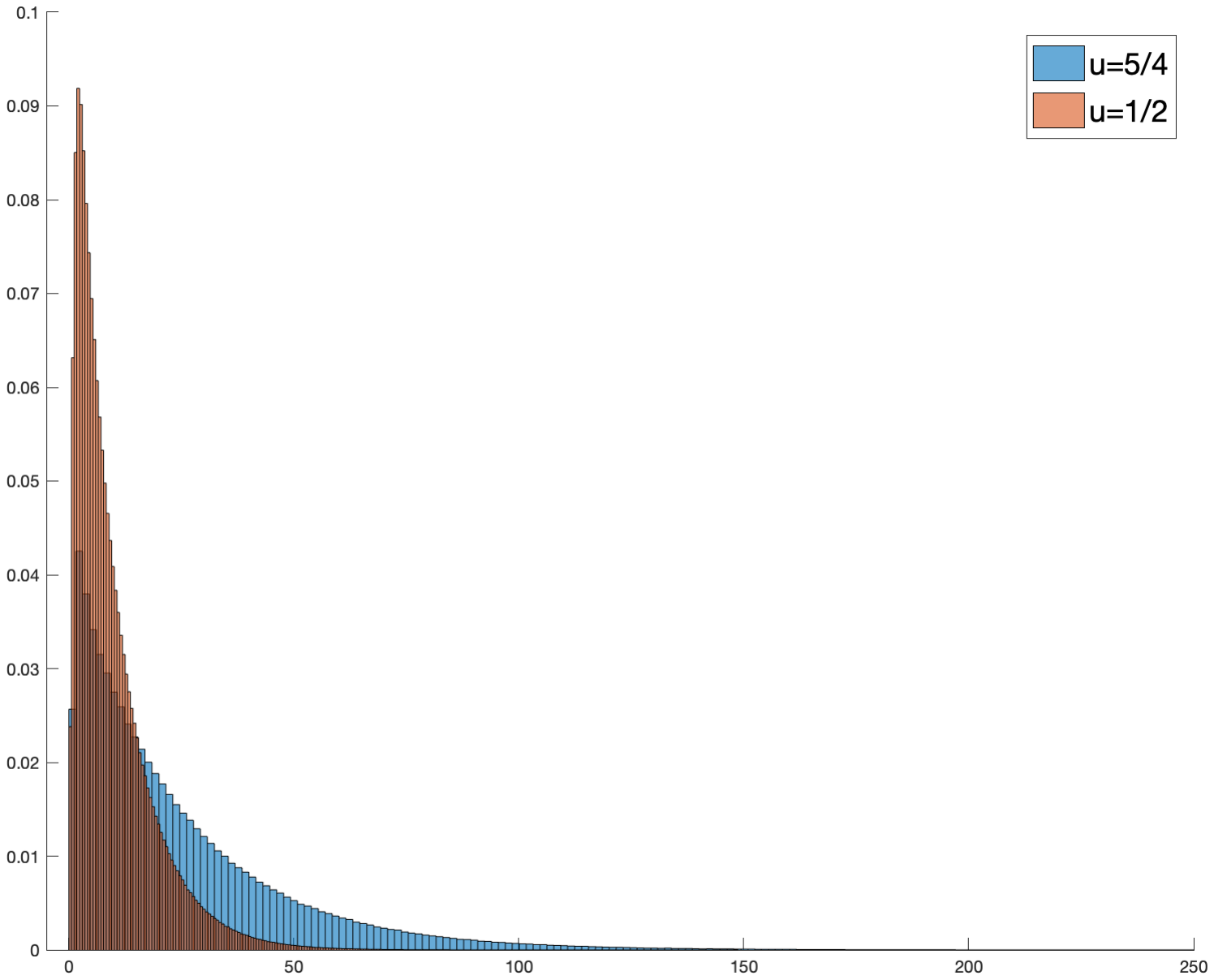}
    \end{minipage}\hfill
\caption{Normalized histograms of the approximated excursion distribution for the Slepian model based on the diffusion covariance, for an $u$ up-crossing {\it (left)} and an $u$ down-crossing {\it (right)}, for $u=1/2,5/4$.}
\label{fig:expdist}
\end{figure}

\begin{table}[h]
\centering
\caption{Simulation estimates of persistency, based on trajectories}
\label{tab:Wafo}
\begin{tabular}{ccccc}
\toprule[0.9pt]
Crossing level, $u$ & Estimate, $\theta_+$ & Confidence interval & Estimate, $\theta_-$ & Confidence interval  \\
\midrule[0.9pt]
$0$ & $0.1885$ & $\pm0.000164$ & $ 0.1883$ & $\pm0.000165$ \\
$1/2$ & $0.2932$ & $\pm0.000271$ & $ 0.1117$ & $\pm0.000104 $\\
$1$ & $0.4295$ & $\pm0.000482$ & $ 0.0598$ & $\pm0.000067 $\\
$5/4$ & $0.5101$ & $\pm0.000661$  & $0.0417$ & $\pm0.000054 $ \\
\bottomrule[0.9pt]
\end{tabular}
\end{table}
It should be noted that Slepian-based estimates of persistency are remarkably close to the estimates based on simulating trajectories. Except for the zero level crossings for which the true persistency coefficient was shown to be $0.1875$ by \cite{PoplavskyiSchehr2018}. For the sake of this example, it can be stated that the Slepian-based IIA provides a good approximation of persistency for Gaussian diffusion processes.

The Slepian IIA not only estimates persistency coefficients but also provides information about the approximated excursion time distribution. In Figure~\ref{fig:expdist}, normalized histograms are presented such that they correspond to the density of the excursion distribution. From Figure~\ref{fig:expdist}, the exponential nature of the tail is evident, and the long tail of the down crossing for higher levels is also clear. It also shows how shifting the level from $1/2$ to $5/4$ alters the tail behavior of the distribution.
Since an arbitrary number of samples can be generated, any characteristics of interest, such as the mean and moments, are approximately obtainable from these samples. This illustrates the usefulness of the Slepian-based IIA.

\section{Conclusion}
\label{sec:conclusion}
\noindent
Using the Slepian-based IIA to approximate the excursion time, distributions have been extended to non-zero crossing levels. By utilizing the crossing behavior of the process of interest and matching the clipped version of this process to a non-stationary switch process, probabilistic valid approximations are obtained. An application to the Gaussian diffusion process is presented, and the results are in line with those obtained from simulations. There are several natural directions for future work on this topic. Perhaps the most natural step is to investigate the asymptotics of the IIA at high crossing levels.

\section{Acknowledgement} 
\noindent
Henrik Bengtsson acknowledges financial support from The Royal Physiographic Society in Lund. Krzysztof Podg\'orski and Henrik Bengtsson also acknowledge the financial support of the Swedish Research Council (VR) Grant DNR: 2020-05168. The authors are also grateful for the clarifying comments from Clément Sire. 

\bibliographystyle{apalike}
\bibliography{References}
\newpage

\appendix
\section*{Appendix: Proofs and auxiliary results}
\noindent
A lemma was needed to derive the expectation of the clipped Slepian process, which is presented below. Recall that the Rayleigh distribution with the parameter $\sigma_{r}^2$  is given by the density $f_R(s)=se^{-s/\sigma_r^2}/\sigma_r^2$, and we shortly write 
$R\sim Rayleigh(\sigma_r)$ when a random variable $R$ has this distribution. 
\begin{lemma}
\label{RN} 
Let $R\sim Rayleigh(\sigma_r)$ and $W\sim  \mathcal{N}(\mu,\sigma_w^2)$ be independent of each other, then 
\begin{align*}
\Xi(\tilde \mu,\tilde \sigma)\stackrel{def}{=}\Prob(R<W)&=\Phi\left(\tilde \mu\right)- \frac{\tilde \sigma}{\sqrt{1 + \tilde \sigma^2}}  \Phi \bigg (\tilde \mu\frac{\tilde \sigma}{\sqrt{1+\tilde \sigma^2}} \bigg) \exp\left( -\frac{\tilde \mu^2}{2(1+\tilde \sigma^2)} \right), 
\end{align*}
where $(\tilde \mu,\tilde \sigma)=(\mu,\sigma_r)/\sigma_w$.
\end{lemma}
\begin{proof}
Without loss of generality, we can assume that $\sigma_r=1$. Under the stated conditions and $F_R$ being the cdf of $R$:
\begin{align*}
\Prob(R<W)&=\Ex(F_R(W))=\Ex\left(\mathbb{I}_{W>0}(1-e^{-{W^2}/{2}})\right) \\
&=P(W>0)-  \frac{1}{\sqrt{2 \pi \sigma_w^2 }}\int_{0}^{\infty}  e^{-\frac{w^2}{2}-\frac{(w-\mu)^2}{2\sigma_w^2}}dw.
\end{align*}
We have 
\begin{align*}
\frac{1}{\sqrt{2 \pi \sigma_w^2}} \int_{0}^{\infty}   e^{- \frac{ (w-\mu)^2 + \sigma_w^2  w^2}{2  \sigma_w^2}}dw 
&=
\frac{1}{\sqrt{2 \pi \sigma_w^2}}  \int_{0}^{\infty}   e^{ {-\frac{w^2-2w\mu/(1 + \sigma_w^2)  + \mu^2/(1 + \sigma_w^2)}{2\sigma_w^2/(1 + \sigma_w^2)}}}dw 
\\
&=\frac{1}{\sqrt{2 \pi \sigma_w^2}} 
\int_{0}^{\infty} e^{-\frac{ \left(w- { \mu}/{(1 + \sigma_w^2)} \right)^2 + \mu^2\left(1/(1+\sigma_w^2)-1/(1+\sigma_w^2)^2\right) }{2\sigma_w^2/(1 + \sigma_w^2) }}dw 
\\
&=\frac{1}{\sqrt{1 + \sigma_w^2}}  e^{-\frac{\mu^2}{2(1+\sigma_w^2)}} \int_{0}^{\infty} \frac{1}{\sqrt{2 \pi \frac{\sigma_w^2}{1+ \sigma_w^2}}}  e^{ -\frac{ \left( w- { \mu}/{(1 + \sigma_w^2)} \right)^2 }{ 2\sigma_w^2/(1 + \sigma_w^2) }}dw 
\\
&=
\frac{1}{\sqrt{1+ \sigma_w^2}}  e^{
-\frac{\mu^2}{2(1+\sigma_w^2)}
} 
\Phi \bigg (\frac{\mu}{\sigma_w}/{\sqrt{1+\sigma_w^2}} \bigg), 
\end{align*}
which yields the result. 
\end{proof}

\begin{corollary}
\label{cor:mu0}
If in the above result $\mu=0$, then 
\begin{align}
\Prob(R<W)=\frac12\left(1-\sqrt{\frac{\sigma_r^2}{\sigma_r^2+\sigma_w^2}}\right).
\end{align}
\end{corollary}

\begin{proof}[Proof of Theorem~\ref{th:expctSlep}]
Let us define a normal variable 
\begin{align}
W=\frac{\sqrt{-r''(0)}}{r'(t)} \left( \Delta(t)+u \frac{r(t)-r(0)}{r(0)} \right).
\end{align}
The mean and variance of this variable are
\begin{align}
\mu=u \sqrt{-r''(0)}\frac{r(t)-r(0)}{r(0) r'(t)},\,\, \sigma_w^2=\frac{r''(0)}{r(0)}\frac{r^2(t)-r^2(0)}{r'^2(t)} -1.
\end{align}
Using Lemma~\ref{RN}, we get 
\begin{align*}
E_u^+(t)&=2\Prob(X_u>u)-1  =2 \Prob(W>R) -1 =2 \Xi(\mu/\sigma_w,1/\sigma_w)-1.
\end{align*}
After some algebra, we get
\begin{align*}
E_u^+(t)&=2\Phi\left(
    \frac{u}{\sqrt{r(0)}} \frac{r(t)-r(0)}{\sqrt{r^2(0)-r^2(t)+\frac{r(0)}{r''(0)}r'^2(t)}} \right)-1
    \\
    &-2\sqrt{-\frac{r(0)}{r''(0)}}\frac{r'(t)}{\sqrt{r^2(0)-r^2(t)}}\exp\left(-\frac{u^2}{2r(0)}\frac{r(0)-r(t)}{r(0)+r(t)} \right)
    \\
    &\times \Phi\left(\frac{-u}{\sqrt{-r''(0)}} \sqrt{\frac{r(0)-r(t)}{r(0)+r(t)}}\frac{r'(t)}{\sqrt{r^2(0)-r^2(t)+\frac{r(0)}{r''(0)}r'^2(t)}}\right).
\end{align*}
\end{proof}

Next, we evaluate the conditional expectation of a clipped Slepian process ( \eqref{eq:clippedSlepianup}) for a given value of the Rayleigh variable $R$.
\begin{lemma}
Let $u>0$ and $D^+_u$ be defined as in  \eqref{eq:clippedSlepianup}.
Then for each $s>0$:
\begin{align*}
\Ex(D^+_u(t) | R=s)&=1 -2\Phi \left( \frac{u\frac{r(0)-r(t)}{r(0)}-s\frac{r'(t)}{\sqrt{-r''(0)}}}{\sqrt{\frac{r^2(0)-r^2(t)}{r(0)}+\frac{r'^2(t)}{r''(0)}}}  \right).
\end{align*}
\end{lemma}
\begin{proof} 
Since conditionally, on $R$, the Slepian process is Gaussian, the derivation is straightforward 
\begin{align*}
E(Y_u(t) | R=s)=&E \big ( I(X_0(t)>u)-I(X_0(t)<u) | R=s \big )  
\\
=&P(ur(t)-Rr'(t)+\Delta(t)>u | R=s) 
\\
&- P(ur(t)-Rr'(t)+\Delta(t)<u | R=s) 
\\
=&1-2P(ur(t)-sr'(t)+\Delta(t)<u | R=s) 
\\
=&1-2P\left(Z <  -\frac{u(r(t)-1)-sr'(t)}{\sqrt{r(0)-r(t)^2-r'(t)^2}} \right) 
\\
=&1-2\Phi \bigg (- \frac{u(r(t)-1)-sr'(t)}{\sqrt{r(0)-r(t)^2-r'(t)^2}}  \bigg ),
\end{align*}
where $Z$ in the above is, as usual, a standard normal variable. 
\end{proof}
In the discussion of the IIA, we used the following result to determine the distribution of the number of switches.
\begin{lemma} 
\label{lemma:stat}
Let $\tilde N(t)$, $t\ge 0$, be the number of switches in $[0,t]$ of a stationary switch process $\tilde D$. Then  
\begin{multline*}
\Prob{\left(\tilde N(t)=k \big | \tilde D(0)= \delta\right)}=\\
=
\begin{cases}
1-F_{A|\delta}(t); k=0 \\
 \left(F_+ \star  F_{-}\right)^{(k/2-1)\star }\star F_{-\delta}\star  F_{A|\delta}(t)-\left(F_+ \star  F_{-}\right)^{(k/2)\star}\star F_{A|\delta}(t);& k>0 \mbox{ is even},
 \\
\left(F_+ \star  F_{-}\right)^{(k-1)/2\star}\star F_{A|\delta}(t)-\left(F_+ \star F_{-}\right)^{(k-1)/2\star}\star F_{-\delta}\star F_{A|\delta}(t);& k \mbox{ is odd},
\end{cases}
\end{multline*}
where $\star$ is used to denote the convolution of probability distributions (i.e., the distribution function of the sum of random variables),   $F_{A|\delta}$  is the cdf of the initial delay conditionally on $\delta$.
\end{lemma}

\begin{proof}
    We consider only the case of conditioning on $\delta=1$ as the opposite case can be obtained by the symmetry argument. 
We first note that for a positive $t$ and a non-negative integer $l$:
\begin{align*}
\Prob{\left(\tilde N(t)=0\big | \delta=1\right)}&=\Prob{\left(A>t|\delta =1\right)}=1- F_{A|\delta =1}(t),
\\
\Prob{\left(\tilde N(t)=1\big | \delta=1\right)}&=\Prob{(A+T^-_1>t\ge A|\delta=1})
\\
&=\Prob{(A+T^-_1>t}|\delta =1)- \Prob{(A+T^-_1>t, A>t}|\delta =1)
\\
&=1-F_{A|\delta=1}\star F_-(t)-\left( 1- F_{A|\delta=1}(t)\right),
\\
\Prob{\left(\tilde N(t)=2l\big | \delta=1\right)}&=\Prob{\left(A+ \sum_{i=1}^{l}{T^+_i}+\sum_{i=1}^{l}{T^-_i}> t \ge A+\sum_{i=1}^{l-1}{T^+_i}+\sum_{i=1}^{l}{T^-_i} \Big |\delta=1\right)}
\\
&=1-F_+^{l\star} \star F_-^{l\star }\star F_{A|\delta=1}(t)-\left(1- F_+^{(l-1)\star}\star F_-^{l\star}\star F_{A|\delta=1}(t)\right)
\\
&= F_+^{(l-1)\star}*F_-^{l\star}\star F_{A|\delta=1}(t)-F_+^{l\star}\star F_-^{l\star}\star F_{A|\delta=1}(t),
\\
\Prob{\left(\tilde N(t)=2l+1\big | \delta=1\right)}&=\Prob{\left( A+\sum_{i=1}^{l}{T^+_i}+\sum_{i=1}^{l+1}{T^-_i} >t \ge  A+\sum_{i=1}^{l}{T^+_i}+\sum_{i=1}^{l}{T^-_i} \right)}
\\
&=F_+^{l\star}\star F_-^{l\star}\star F_{A|\delta=1}(t)-F_+^{l\star} \star F_-^{(l+1)\star}\star F_{A|\delta=1}(t).
\end{align*}
\end{proof}
For the next proposition, we need the Laplace transform of $f_{A|\delta}$ and $F_{A|\delta}$. 
\begin{align*}
    \Psi_{A\vert \delta}(s)&=\mathcal{L}(f_{A|\delta})(s)=\frac{1-s\mathcal L (F_\delta)(s)}{s\mu_\delta}=\frac{1-\Psi_\delta(s)}{s\mu_\delta}
    \\
    \mathcal{L}(F_{A|\delta})(s)&=\frac{\mathcal{L}(f_{A|\delta})(s)}{s}=\frac{1-\Psi_\delta(s)}{s^2\mu_\delta}
\end{align*}

\begin{proposition}
\label{prop:exp}
For a stationary switch process, the average number of switches in $[0,t]$, $ N_{<}(t)=\Ex(\tilde N(t)|\delta =-1)$, $ N_{>}(t)=\Ex(\tilde N(t)|\delta =1)$ have Laplace transforms 
\begin{align*}
    \mathcal L{ (N_{>})}(s)&=\frac{1 }{s^2\mu_+} \frac{ \left( 1 +\Psi_{-}(s)\right) \left(1-\Psi_{+}(s)\right)}{1-\Psi_+(s)\Psi_-(s)},
    \\
     \mathcal L{ (N_{<})}(s)&=\frac{1}{s^2\mu_-}\frac{ \left( 1 + \Psi_{+}(s)\right) \left(1-\Psi_{-}(s)\right)}{1-\Psi_+(s)\Psi_-(s)}.
\end{align*}
\end{proposition}
\begin{proof}
The result follows from 
\begin{align*}
    \mathcal L{ (N_{>})}(s) &=\sum_{k=1}^\infty k\cdot  \mathcal L \left(\Prob\left(\tilde N(\cdot)=k|\delta=1\right) \right)(s)
    \\
    &=\sum_{l=0}^\infty (2l+1) \left( \mathcal L\left( \left(F_+ \star  F_{-}\right)^{l\star}\star F_{A|\delta=1}\right)(s) -\mathcal L\left(\left(F_+ \star F_{-}\right)^{l\star}\star F_{-}\star F_{A|\delta=1}\right)(s)  \right)
    \\
    &+ \sum_{l=1}^\infty 2l \left(\mathcal L\left( \left(F_+ \star  F_{-}\right)^{(l-1)\star }\star F_{-}\star  F_{A|\delta=1}\right)(s) -\mathcal L\left( \left(F_+ \star  F_{-}\right)^{l\star}\star F_{A|\delta=1}\right) (s) \right)
    \\
    &=\frac{\Psi_{A|\delta=1}}{s} \left( \left( 1 -\Psi_{-}\right) \sum_{l=0}^\infty (2l+1) \left( \Psi_+ \Psi_{-}\right)^{l} + \left( 1 - \Psi_{+}\right)\sum_{l=1}^\infty 2l \left(\Psi_+  \Psi_{-}\right)^{(l-1) }\Psi_{-}  \right)
    \\
    &=\frac{1 - \Psi_{+}}{s^2\mu_+} \left( \frac{1 - \Psi_{-}}{1-\Psi_+\Psi_-}+ 2\Psi_-\sum_{l=1}^\infty l \left(\Psi_+ \Psi_{-}\right)^{l-1} \left( \left( 1 - \Psi_{-}\right)\Psi_+ + 1 -  \Psi_{+} \right)\right)
 \\
  &=\frac{ 1 - \Psi_{+}}{s^2\mu_+} \left( \frac{ 1 - \Psi_{-}}{1-\Psi_+\Psi_-}+  \frac{2\Psi_- \left(1- \Psi_{-}\Psi_+ \right)}{\left(1-\Psi_+\Psi_-\right)^2}  \right)
  =\frac{ 1 -  \Psi_{+}(s)}{s^2\mu_+} \frac{ 1 +  \Psi_{-}(s)}{1-\Psi_+(s)\Psi_-(s)}.
\end{align*}
\end{proof}





\end{document}